\documentclass[12pt]{amsart}

\usepackage{amsthm,amssymb,amsmath,amstext,amsfonts}
\usepackage{enumitem,mathtools,pgfplots,pgfmath,subcaption}
\usepackage{caption}
\usepackage{tikz,subcaption}
\tikzset{>=latex}
\usepackage{tipa}
\usepackage{colonequals}
\usepackage{dirtytalk}
\usepackage{gensymb}
\usepackage{multicol}
\usetikzlibrary{hobby,arrows,shapes,automata,backgrounds,decorations,petri,positioning,arrows.meta}
\usetikzlibrary{decorations.pathreplacing,angles,quotes}
\pgfdeclarelayer{background}
\pgfsetlayers{background,main}
\usepackage[hyphens]{url}
\usepackage{float}
\usepackage{bbm}
\pgfdeclarelayer{background}
\pgfsetlayers{background,main}

\usepackage[margin=1in,letterpaper,portrait]{geometry}

\usepackage[latin1]{inputenc}
\usepackage{url}
\usepackage{chngcntr}
\usepackage{apptools}
\usepackage{amsaddr}
\AtAppendix{\counterwithin{lemma}{section}}

\theoremstyle{plain}
\theoremstyle{definition}
\newtheorem{theorem}{Theorem}[section]
\newtheorem*{theorem-nonum}{Theorem}
\newtheorem{conjecture}[theorem]{Conjecture}

\newtheorem{lemma}[theorem]{Lemma}
\newtheorem{definition}[theorem]{Definition}
\newtheorem*{definition-nonum}{Definition}
\newtheorem{question}[theorem]{Question}
\newtheorem{example}[theorem]{Example}

\newtheorem{corollary}[theorem]{Corollary}
\newcommand*{\Scale}[2][4]{\scalebox{#1}{$#2$}}

\DeclareMathAlphabet{\mathpzc}{OT1}{pzc}{m}{it}
\DeclareMathOperator{\cost}{cost}
\DeclareMathOperator{\lps}{lps}
\DeclareMathOperator{\red}{red}

\begin{document}

\title{Sorting by shuffling methods and a queue}

\author{Stoyan Dimitrov}
\address{University of Illinois at Chicago}
\email{sdimit6@uic.edu}



\begin{abstract}
We study sorting by queues that can rearrange their content by applying permutations from a predefined set. These new sorting devices are called \emph{shuffle queues} and we investigate those of them corresponding to sets of permutations defining some well-known shuffling methods. If $\mathbb{Q}_{\Sigma}$ is the shuffle queue corresponding to the shuffling method $\Sigma$, then we find a number of surprising results related to two natural variations of shuffle queues denoted by $\mathbb{Q}_{\Sigma}^{\prime}$ and $\mathbb{Q}_{\Sigma}^{\textsf{pop}}$. These require the entire content of the device to be unloaded after a permutation is applied or unloaded by each pop operation, respectively.

First, we show that sorting by a deque is equivalent to sorting by a shuffle queue that can reverse its content. Next, we focus on sorting by \emph{cuts}. We prove that the set of permutations that one can sort by using $\mathbb{Q}_{\text{cuts}}^{\prime}$ is the set of the $321$-avoiding separable permutations. We give lower and upper bounds to the maximum number of times the device must be used to sort a permutation. Furthermore, we give a formula for the number of $n$-permutations, $p_{n}(\mathbb{Q}_{\Sigma}^{\prime})$, that one can sort by using $\mathbb{Q}_{\Sigma}^{\prime}$, for any shuffling method $\Sigma$, corresponding to a set of irreducible permutations.

We also show that $p_{n}(\mathbb{Q}_{\Sigma}^{\textsf{pop}})$ is given by the odd indexed Fibonacci numbers $F_{2n-1}$, for any shuffling method $\Sigma$ having a specific \say{back-front} property. The rest of the work is dedicated to a surprising conjecture inspired by Diaconis and Graham, which states that one can sort the same number of permutations of any given size by using the devices $\mathbb{Q}_{\text{In-sh}}^{\textsf{pop}}$ and $\mathbb{Q}_{\text{Monge}}^{\textsf{pop}}$, corresponding to the popular \emph{In-shuffle} and \emph{Monge} shuffling methods.
\end{abstract}

\maketitle

\section{Introduction and definitions}
If we have a device that can rearrange the elements of a given input permutation according to certain rules, then a natural question is: \say{Which permutations of $1,2, \cdots ,n$ can be sorted when we use this device?} The first to consider a question with this formulation was Tarjan \cite{tarjan} who, like others, was inspired by one chapter in the seminal book of Donald Knuth \cite[Chapter 2.2.1]{knuth}. 

\begin{figure}[h!]
\centering
\begin{subfigure}{0.3\textwidth}
\centering
           \begin{tikzpicture}[scale=0.22]
           \centering
\draw[black, thick] (-1.5,7) -- (-1.5,0); \draw[black, thick] (1.5,0) -- (1.5,7);
\node (inp) at (8,3.5) {Input};
\node (out) at (-8,3.5) {Output};
\node (devInp2) at (0.75,6) {};
\node (devOut1) at (-0.75,6) {};

\draw [->] (devOut1) to [out=90, in=90] node [above=0.8mm] (TextNode1) {$O$} (out);
\draw [->] (inp) to [out=90, in=90] node [above=0.66mm] (TextNode1) {$I$} (devInp2);
\end{tikzpicture}
\vspace*{3.5mm}
\caption{stack}
\end{subfigure}
\begin{subfigure}{0.3\textwidth}
\centering
         \begin{tikzpicture}[scale=0.22]
         \centering
\draw[black, thick] (-1.5,7) -- (-1.5,0); \draw[black, thick] (1.5,0) -- (1.5,7);
\node (inp) at (8,3.5) {Input};
\node (out) at (-8,3.5) {Output};
\node (devOut1) at (-0.75,6) {};
\node (devInp1) at (0.75,1) {};

\draw [->] (devOut1) to [out=90, in=90] node [above=0.8mm] (TextNode1) {$O$} (out);

\draw [->] (inp) to [out=270, in=270] node [below=0.66mm] (TextNode1) {$I$} (devInp1);
\end{tikzpicture}
\caption{queue}
\end{subfigure}
\begin{subfigure}{0.3\textwidth}
\centering
         \begin{tikzpicture}[scale=0.22]
         \centering
\draw[black, thick] (-1.5,7) -- (-1.5,0); \draw[black, thick] (1.5,0) -- (1.5,7);
\node (inp) at (8,3.5) {Input};
\node (out) at (-8,3.5) {Output};
\node (devInp2) at (0.75,6) {};
\node (devOut1) at (-0.75,6) {};
\node (devInp1) at (0.75,1) {};
\node (devOut2) at (-0.75,1) {};

\draw [->] (devOut1) to [out=90, in=90] node [above=0.8mm] (TextNode1) {$O$} (out);
\draw [->] (inp) to [out=90, in=90] node [above=0.66mm] (TextNode1) {$\overline{I}$} (devInp2);

\draw [->] (devOut2) to [out=270, in=270] node [below=0.8mm] (TextNode1) {$\overline{O}$} (out);
\draw [->] (inp) to [out=270, in=270] node [below=0.66mm] (TextNode1) {$I$} (devInp1);
\end{tikzpicture}
\caption{deque}
\label{fig:3figsD}
\end{subfigure}
     \caption{The input and output operations on stack, queue and deque}
     \label{fig:3devices}
 \end{figure}
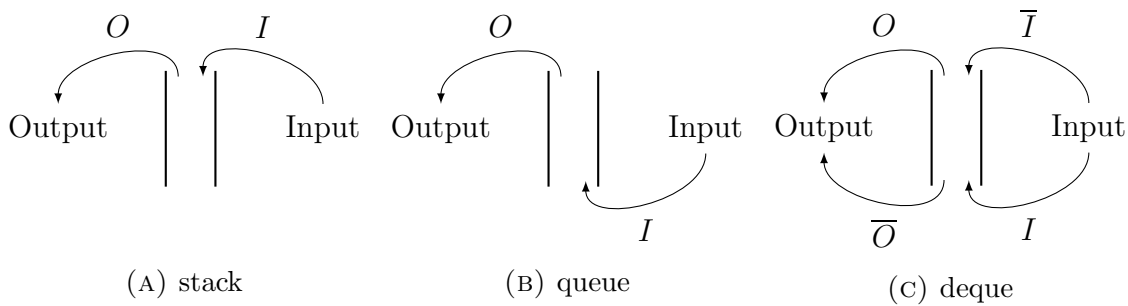

Knuth considered the classical data structures \emph{stack}, \emph{queue} and \emph{deque} (double-ended queue) which are shown at Figure \ref{fig:3devices}. He asked which permutations can be obtained by using each of these devices, if we begin with the identity permutation $12\cdots n$. The two questions correspond to two equivalent viewpoints since a permutation $\pi$ can be obtained from the identity by applying a given sequence of operations, if and only if $\pi^{-1}$ is sorted by the same sequence of operations.

Below is a brief description of the three devices. All of them are linear lists which are used frequently in programming to store and access data. For each device, we have input operations (also called \emph{push} operations), which insert an element from the input to the device and output operations (also called \emph{pop} operations), which move an element from the device to the output:
\begin{itemize}
    \item stack ($\mathbbm{St}$): the input operations $I$ and the output operations $O$ are made at one end of the list.
    \item queue ($\mathbbm{Q}$): the input operations $I$ are made at one end of the list and the output operations $O$ are made at the other end of the list.
    \item deque ($\mathbbm{Deq}$): two kinds of input operations ($I$ and $\overline{I}$) exist, as well as two kinds of output operations ($O$ and $\overline{O}$). The two pairs of input and output operations are made at the two opposite ends of the list, as shown at Figure \ref{fig:3figsD}.
\end{itemize}

The question of Knuth led to the development of permutation pattern research. The most cited result of his work is the fact that the permutations that can be sorted with a stack are the $231$-avoiding permutations. A great number of subsequent articles investigated sorting by different variations of a stack or networks of stacks. Some examples are \emph{pop-stacks} \cite{avis}, \emph{stacks in parallel} \cite{itai,tarjan}, \emph{stacks in series} \cite{west,Z} and \emph{stacks of bounded size} \cite{boundedStacks}. 

Sorting by a deque and its variations has also been a subject of serious research interest. It was proved by Pratt \cite{pratt} that the deque sortable permutations are characterized by avoiding a certain infinite set of permutations and the enumeration of this set is still an open problem. Knuth himself \cite[Chapter 2.2.1]{knuth} considered \emph{input-restricted} and \emph{output-restricted} deques and determined the sets of obtainable (respectively, sortable) permutations by them and his approach introduced the important \emph{kernel method}. An interesting recent result of Price \cite{price} states that the permutations sortable by a deque and two stacks in parallel share the same growth rate. 

The only permutation that can be sorted by a queue is the identity permutation and this, as Knuth writes, follows trivially \say{by the nature of the queue.} One of the few articles discussing sorting by modification of the queue is \cite{doyle}, where the author looks at a queue that is capable of doing direct transfers of elements from the input to the output. The permutations that can be sorted with this device are the $321$-avoiding permutations. Albert et. al \cite{Cmachines} consider a more general type of devices, called $\mathcal{C}$-machines, that can perform the same direct transfers of elements. For more background on sorting devices, we refer to the surveys \cite{sortingSurvey} and \cite[Chapter 4]{50years}, as well as to the books \cite[Chapter 8]{Bona} and \cite[Chapter 2]{KitBook}.

A completely different, yet connected, line of research investigates \emph{shuffling methods} for a given deck (to not be confused with \say{deque}) of cards or a given permutation. A shuffling method is a procedure that will lead to a uniformly shuffled deck after applying it multiple times. This procedure is usually comprised of the following two steps: choose a permutation out of a given set and then apply it over the deck. Thus any shuffling method has a set of permutations associated with it. Diaconis, Fulman and Holmes \cite[Section 2.3]{diaconisMain} give an overview of the previous work related to shuffling. The mathematics of shuffling uses tools related to mixing times \cite{diaconisThingsUp}, representation theory \cite{odl} and quasi-symmetric polynomials \cite{stanley}. 

In this work, we relate the areas of sorting devices and shuffling methods by considering sorting by special type of queues, called \emph{shuffle queues}, which can rearrange their content by applying permutations in a given collection over it. We will call any such collection of permutations a \emph{shuffling method} and we will focus on collections associated with some methods that are popular in the literature. Shuffle queues are very similar to the \emph{permuting machines} introduced in a paper of Albert et al. \cite{permutingMachines}. However, our settings are more general, since the permuting machines have to satisfy one important property which does not necessary hold in the case of shuffle queues. The two concepts will be equivalent if we require the set of permutations associated with each shuffle queue to be closed under pattern containment. 

Except for the few studies mentioned above, not many previous works investigate sorting by modifications of a queue. Shuffle queues are a natural such modification since a sorting device is a machine whose sole function is to re-order its input data. These new devices lead to some surprising enumerative results and raise interesting combinatorial questions. More motivational points are described in Section \ref{sec:motivation}. 

\subsection{Notation} 
\label{sec:defs}
 The set of consecutive integers $\{i, i+1, \dots , j\}$ will be denoted by $[i,j]$. A permutation of size $n$ is a bijective map from $[n] \coloneqq [1,n]$ to itself. When referring to permutations, we will use their one-line representation. The set of all permutations of size $n$ will be denoted by $S_{n}$. If $\lambda$ is a sequence of distinct numbers, the reduction of $\lambda$, denoted $\red(\lambda)$, is the permutation obtained from $\lambda$ by replacing its $i$-th smallest entry by $i$. For example, we have $\red(4968) = 1423$. A permutation $\pi$ contains a permutation, or a \emph{pattern}, $\sigma$ if there is a subsequence $\lambda$ of $\pi$ such that $\red(\lambda) = \sigma$. If $\pi$ does not contain $\sigma$, then $\pi$ \emph{avoids} $\sigma$. We denote the set of permutations of size $n$ that avoid all the patterns in a set $X$ by $Av_{n}(X)$, with $av_{n}(X)\coloneqq |Av_{n}(X)|$. If $x$ and $y$ are sequences of integers, then we will write $x>y$ (respectively, $x<y$) if each element of $x$ is greater (respectively, less) than each element of $y$. A \emph{segment} of a permutation $\pi = \pi_{1}\cdots\pi_{n}$ will be a subsequence $\pi_{a}\pi_{a+1}\cdots\pi_{b}$ of consecutive elements of $\pi$, for some $1\leq a < b \leq n$ and it will be denoted by $[a,b]$ when $\pi$ is inferred from the context. A \emph{permutation class} $\mathcal{C}$ is a set of permutations, such that if $\pi\in \mathcal{C}$ and $\pi$ contains $\sigma$, then $\sigma\in \mathcal{C}$. Other standard definitions related to permutation patterns that will be used can be found in \cite{bevan}. 
 
 The empty sequence will be denoted by $\varepsilon$. For a sequence of distinct numbers $s$, denote by $Im(s)$ the set of elements of $s$ and let $Im(s_{1}, \dots , s_{r}) = \bigcup_{k=1}^{r} Im(s_{k})$. Consider the set of triples of sequences that partition $[n]$, 
\begin{equation*}
    \{(s_{1},s_{2},s_{3})\mid Im(s_{i})\cap Im(s_{j}) = \emptyset , Im(s_{1},s_{2},s_{3}) = [n]\}.
\end{equation*} 
We will call the elements of this set, for any positive integer $n$, \emph{configurations}.

A sorting device $\mathbb{D}$ is a tool that transforms a given input permutation $\pi$ by following a particular algorithm which could be deterministic or non-deterministic. The result is an output permutation $\pi'$. During the execution of the algorithm, every device $\mathbb{D}$ has a given configuration $(s_{inp},s_{dev},s_{out})$, comprised of three sequences (strings) corresponding to the current string in the input, in the device, and in the output, respectively. The initial configuration is $(\pi, \varepsilon , \varepsilon)$ and the final configuration is $(\varepsilon , \varepsilon , \pi')$. Denote by $\mathbb{D}(\pi)$ the set of possible output permutations, when using a device $\mathbb{D}$ on input $\pi$. If $id_{n}$ denotes the identity permutation of size $n$, then let $S_{n}(\mathbb{D}) \coloneqq \{\pi \mid \pi\in S_{n}, id_{n}\in \mathbb{D}(\pi)\}$ be the set of the permutations sortable by $\mathbb{D}$. Furthermore, let $p_{n}(\mathbb{D}) \coloneqq |S_{n}(\mathbb{D})|$.

In this paper, a shuffling method $\Sigma$ is defined by a family of sets of permutations $\{\Pi_{\Sigma}^{n}\subseteq S_{n}\mid n = 2,\dots\}$ that one can apply over the content of a sorting device, when using the method. Note that $\Pi_{\Sigma}^{n}$ contains permutations of size $n$, for every $n\geq 2$. We will also assume that $id_{n}\notin \Pi_{\Sigma}^{n}$, for every $n\geq 2$. We will refer to $\{\Pi_{\Sigma}^{n}\}_{n=2}^{\infty}$ as the \emph{permutation family} of the method $\Sigma$. We will also use the notations $\Pi( \Sigma )\coloneqq \bigcup_{n=2}^{\infty}\Pi_{\Sigma}^{n}$ and $(\Pi_{\Sigma}^{k})^{-1} \coloneqq \{\sigma^{-1}\mid \sigma\in \Pi_{\Sigma}^{k}\}$. An example of a shuffling method is shuffling by cuts, which is defined below and studied in Sections \ref{sec:cuts} and \ref{sec:higherCost}.
\begin{definition} 
\label{def:cuts}
The shuffling method $cuts$.
    \begin{equation}
    \label{eq:defFcuts}
        \forall n\geq 2: \Pi_{\text{cuts}}^{n} \coloneqq \{ k(k+1)\cdots n12\cdots (k-1)\mid k\in [2,n]\}.
    \end{equation}
\end{definition}

In the existing literature, a shuffling method transforms a given input permutation by multiplying it by another permutation, according to a given distribution. For example, when one uses shuffling by \emph{cuts} over an input of size $n$, one picks a permutation in the set $\{ k(k+1)\cdots n12\cdots (k-1)\mid k\in [2,n]\}$ according to uniform distribution and applies it over the input. In the present paper, we will be able to choose the permutation that can be applied to the content deterministically and thus our definition of a shuffling method does not involve a distribution.

For a given shuffling method $\Sigma$, we consider a non-deterministic sorting device $\mathbb{Q}_{\Sigma}$ for which at any given step one can apply up to three possible operations over the current configuration $s=(s_{inp},s_{dev},s_{out})$. Denote the next configuration by $\overline{s}$. The three operations are described below.

\begin{enumerate}[label=\emph{\arabic*}.]
    \item \textbf{Push}
    
    Move the first element $x$ of the input $s_{inp} = xs'_{inp}$ to the content of the device. We get $\overline{s} = (s'_{inp},s_{dev}x,s_{out})$. One can apply this operation only if $s_{inp}\neq \varepsilon$.\\
    
    \item \textbf{Pop}
    
    Move the first element $y$ of the content of the device $s_{dev} = ys'_{dev}$ to the output. We get $\overline{s} = (s_{inp},s'_{dev},s_{out}y)$. One can apply this operation only if $s_{dev}\neq \varepsilon$.\\
    
    \item \textbf{Shuffle} 
    
    Choose a permutation $\sigma\in\Pi_{\Sigma}^{m}$ and apply it over the content of the device $s_{dev}$, where $|s_{dev}|=m$. We get $\overline{s} = (s_{inp}, \sigma s_{dev},s_{out})$. One can apply this operation only if $m\geq 2$ and if the last operation that has been applied is not a shuffle operation.
\end{enumerate}

Note that the device $\mathbb{Q}_{\Sigma}$ functions as a queue since it can receive entries on one of its ends and release entries on the other end. In addition, the content of this queue can be shuffled and thus we will call it a \emph{shuffle queue}. When a certain permutation is chosen to be applied on a shuffle operation, we will say that the shuffle operation is associated with this permutation. Also, note that the restriction to not have two consecutive shuffle operations is reasonable since if one allows applying multiple consecutive shuffle operations for a shuffle queue $\mathbb{Q}_{\Sigma}$, then sorting by this queue would be equivalent to sorting by a queue $\mathbb{Q}_{\Sigma'}$, for which two consecutive shuffle operations are not allowed. Here, $\Sigma'$ would be the shuffling method for which $\Pi_{\Sigma'}^{n} = \langle \Pi_{\Sigma}^{n} \rangle$, for every $n\geq 2$, where $\langle T \rangle$ denotes the subgroup generated by the set $T$. 

Our work focuses on two natural variations of the devices $\mathbb{Q}_{\Sigma}$ that will be called shuffle queues of type $(i)$ and type $(ii)$. They are obtained after imposing two additional restrictions:

\begin{enumerate}[label=\emph{(\roman*})]
    \item The entire content of the device must be unloaded after each shuffle. 
    
    Denote the corresponding sorting device by $\mathbb{Q}_{\Sigma}^{\prime}$.
    
    \item The entire content of the device must be unloaded by each pop operation. 
    
    Denote the corresponding sorting device by $\mathbb{Q}_{\Sigma}^{\textsf{pop}}$. This is the pop-version of the device $\mathbb{Q}_{\Sigma}$ in analogy to the pop version of the stack-sorting device first considered by Avis and Newborn in \cite{avis}. We will also call them \emph{pop shuffle queues}.
\end{enumerate}

Consider the device of type $(i)$, $\mathbb{Q}_{\text{cuts}}^{\prime}$. Example \ref{ex:cutIter} shows one possible sequence of configurations for $\mathbb{Q}_{\text{cuts}}^{\prime}$ and the corresponding operations when sorting the permutation $213564$ with it. Each configuration is written in the column form $\begin{pmatrix}
s_{inp}\\
s_{dev}\\
s_{out}\\
\end{pmatrix}$. In general, if $\mathbb{D}$ is a sorting device and $\pi'\in\mathbb{D}(\pi)$, then any sequence of configurations for $\mathbb{D}$ that begins with $(\pi, \varepsilon , \varepsilon)$ and ends with $(\varepsilon , \varepsilon , \pi')$, together with the list of corresponding operations, will be called an \emph{iteration} of $\mathbb{D}$ over the input $\pi$.

\begin{example}
\label{ex:cutIter}
    Iteration of $\mathbb{Q}_{\text{cuts}}^{\prime}$ over $213645$.
\begin{equation*}
\begin{split}
    \begin{pmatrix}
213645\\
\varepsilon\\
\varepsilon\\
\end{pmatrix}
\xrightarrow{\textsf{push}}
\begin{pmatrix}
13645\\
2\\
\varepsilon\\
\end{pmatrix}
\xrightarrow{\textsf{push}}
\begin{pmatrix}
3645\\
21\\
\varepsilon\\
\end{pmatrix}
\xrightarrow[\text{+unload}]{\substack{\textsf{shuffle}\\(cut)}}
\begin{pmatrix}
3645\\
\varepsilon\\
12\\
\end{pmatrix} 
\xrightarrow{\textsf{push}}
\begin{pmatrix}
645\\
3\\
12\\
\end{pmatrix} \\[20pt]
\xrightarrow{\textsf{pop}}
\begin{pmatrix}
645\\
\varepsilon\\
123\\
\end{pmatrix}
\xrightarrow{\textsf{push}}
\begin{pmatrix}
45\\
6\\
123\\
\end{pmatrix}
\xrightarrow{\textsf{push}}
\begin{pmatrix}
5\\
64\\
123\\
\end{pmatrix}
\xrightarrow{\textsf{push}}
\begin{pmatrix}
\varepsilon\\
645\\
123\\
\end{pmatrix}
\xrightarrow[\text{+unload}]{\substack{\textsf{shuffle}\\(cut)}}
\begin{pmatrix}
\varepsilon\\
\varepsilon\\
123456\\
\end{pmatrix}\\[5pt]
\end{split}
\end{equation*}
\end{example}

This device requires that we unload the entire content of the device after each shuffle operation. Also, note that one can choose to apply multiple different cuts on each shuffle operation. 
Consider the device of type $(ii)$, $\mathbb{Q}_{\text{cuts}}^{\textsf{pop}}$. Below is shown one possible iteration of $\mathbb{Q}_{\text{cuts}}^{\textsf{pop}}$. 

\begin{example}
\label{ex:cutIterPop}
    Iteration of $\mathbb{Q}_{\text{cuts}}^{\textsf{pop}}$ over $41325$.
\begin{equation*}
\begin{split}
    \begin{pmatrix}
41325\\
\varepsilon\\
\varepsilon\\
\end{pmatrix}
\xrightarrow{\textsf{push}}
\begin{pmatrix}
1325\\
4\\
\varepsilon\\
\end{pmatrix}
\xrightarrow{\textsf{push}}
\begin{pmatrix}
325\\
41\\
\varepsilon\\
\end{pmatrix}
\xrightarrow{\textsf{push}}
\begin{pmatrix}
25\\
413\\
\varepsilon\\
\end{pmatrix} 
\xrightarrow{\substack{\textsf{shuffle}\\(cut)}}
\begin{pmatrix}
25\\
341\\
\varepsilon\\
\end{pmatrix} 
\xrightarrow{\textsf{push}}
\begin{pmatrix}
5\\
3412\\
\varepsilon\\
\end{pmatrix}\\[20pt]
\xrightarrow{\substack{\textsf{shuffle}\\(cut)}}
\begin{pmatrix}
5\\
1234\\
\varepsilon\\
\end{pmatrix}
\xrightarrow{\substack{\textsf{pop}\\(unload)}}
\begin{pmatrix}
5\\
\varepsilon\\
1234\\
\end{pmatrix}
\xrightarrow{\textsf{push}}
\begin{pmatrix}
\varepsilon\\
5\\
1234\\
\end{pmatrix}
\xrightarrow{\substack{\textsf{pop}\\(unload)}}
\begin{pmatrix}
\varepsilon\\
\varepsilon\\
12345\\
\end{pmatrix}\\[5pt]
\end{split}
\end{equation*}
\end{example}

The device in Example \ref{ex:cutIterPop} requires that we unload the entire content of the device by each pop operation, but we do not have to do that after a shuffle operation.

\subsection{Motivation}
\label{sec:motivation}

Here, we describe some additional motivation to consider sorting by shuffle queues, as well as their variations of types $(i)$ and $(ii)$. We also motivate the investigation of sorting by cuts, which is a main focus of the present work.

Sorting by a deque is equivalent to sorting by a simple shuffle queue (see Section \ref{sec:deque}). Perhaps, one could find shuffle queues that mirror sorting by other popular devices. This would give new perspectives and might help solving certain problems related to these devices. In addition, sorting by $\mathbb{Q}_{\text{cuts}}$ has a simple interpretation in terms of railway switching networks, which was the way used by Knuth in \cite{knuth} to illustrate sorting by stack, queue and deque. Add a circular railroad extension connecting the beginning and the end of a railroad queue, as on Figure \ref{fig:curcularExt} below.

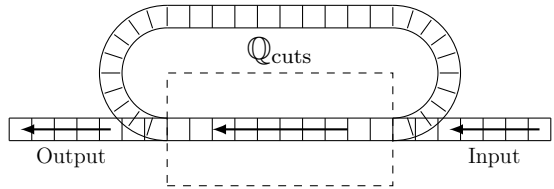
\begin{figure}[H]
  \centering
   \begin{tikzpicture}[every node/.style={inner sep=0pt},scale=0.3]
        \draw (-12,1) -- (12,1);
        \draw (-12,0) -- (12,0);
        
        \foreach \i in {0,...,24}
        {
        \pgfmathtruncatemacro{\y}{-12+\i };
        \draw (\y,0) -- (\y,1);
        }
    \draw (-5,5) -- (5,5);
    \draw (-5,6) -- (5,6);
     \foreach \i in {0,...,8}
        {
        \pgfmathtruncatemacro{\y}{-4+\i };
        \draw (\y,5) -- (\y,6);
        }
    
 \draw (-5,5) arc (90:270:2) 
 node[pos=0] (a0) {} 
 node[pos=0.1] (a1) {} 
 node[pos=0.2] (a2) {}
 node[pos=0.3] (a3) {}
 node[pos=0.4] (a4) {}
 node[pos=0.5] (a5) {} 
 node[pos=0.6] (a6) {}
 node[pos=0.7] (a7) {} 
 node[pos=0.8] (a8) {}
 node[pos=0.9] (a9) {};

 \draw (-5,6) arc (90:270:3)
 node[pos=0] (b0) {} 
 node[pos=0.1] (b1) {} 
 node[pos=0.2] (b2) {}
 node[pos=0.3] (b3) {}
 node[pos=0.4] (b4) {}
 node[pos=0.5] (b5) {}
 node[pos=0.6] (b6) {}
 node[pos=0.7] (b7) {}
 node[pos=0.8] (b8) {}
 node[pos=0.9] (b9) {};

\draw (a0) -- (b0);
\draw (a1) -- (b1);
\draw (a2) -- (b2);
\draw (a3) -- (b3);
\draw (a4) -- (b4);
\draw (a5) -- (b5);
\draw (a6) -- (b6);
\draw (a7) -- (b7);
\draw (a8) -- (b8);
\draw (a9) -- (b9);

 \draw (5,5) arc (90:-90:2)
 node[pos=0] (c0) {} 
 node[pos=0.1] (c1) {} 
 node[pos=0.2] (c2) {}
 node[pos=0.3] (c3) {}
 node[pos=0.4] (c4) {}
 node[pos=0.5] (c5) {} 
 node[pos=0.6] (c6) {}
 node[pos=0.7] (c7) {} 
 node[pos=0.8] (c8) {}
 node[pos=0.9] (c9) {};
 
 \draw (5,6) arc (90:-90:3)
 node[pos=0] (d0) {} 
 node[pos=0.1] (d1) {} 
 node[pos=0.2] (d2) {}
 node[pos=0.3] (d3) {}
 node[pos=0.4] (d4) {}
 node[pos=0.5] (d5) {}
 node[pos=0.6] (d6) {}
 node[pos=0.7] (d7) {}
 node[pos=0.8] (d8) {}
 node[pos=0.9] (d9) {};

\draw (c0) -- (d0);
\draw (c1) -- (d1);
\draw (c2) -- (d2);
\draw (c3) -- (d3);
\draw (c4) -- (d4);
\draw (c5) -- (d5);
\draw (c6) -- (d6);
\draw (c7) -- (d7);
\draw (c8) -- (d8);
\draw (c9) -- (d9);

  \draw[dashed] (-5,-2) rectangle (5,3);
  \node (dev) at (0,4) {$\mathbb{Q}_{\text{cuts}}$};
  
  \draw [<-,thick] (-3,0.5) -- (3,0.5);
  
  \draw [<-,thick] (7.5,0.5) -- (11.5,0.5);
  \node[scale=0.7] (inp) at (9.5,-0.75) {Input};
  \draw [<-,thick] (-11.5,0.5) -- (-7.5,0.5);
  \node[scale=0.7] (inp) at (-9.25,-0.75) {Output};
    \end{tikzpicture}

  \caption{The shuffle queue $\mathbb{Q}_{\text{cuts}}$ represented as a railway switching network.}
  \label{fig:curcularExt}
\end{figure}
Suppose that a railroad car cannot enter or leave the queue (no pushes or pops are allowed), while there is a car in the extension. Thus we have a queue that can move a group of consecutive elements from its beginning to its end. This is exactly what one can do by cuts.

It is not difficult to show that one can sort every permutation using $\mathbb{Q}_{\text{cuts}}$ (Corollary \ref{cor:cutsPop} gives even a stronger statement). Thus, it is reasonable to ask which permutations can be sorted by cuts and by other methods if we consider the two natural restrictions defining shuffle queues of types $(i)$ and $(ii)$, namely, to unload the content after each shuffle or by each pop, respectively. Sorting by the shuffle queue of type $(i)$, $\mathbb{Q}_{\text{cuts}}^{\prime}$, corresponds to sorting by the same railway switching network shown at Figure \ref{fig:curcularExt}, with the additional requirement that we have to unload the queue after each use of the extension.

 $\mathbb{Q}_{\text{cuts}}^{\prime}$ is a non-deterministic device and we show that by using this device one can sort a subset of the separable permutations defined at the beginning of Section \ref{sec:cuts}. Therefore, there exists a deterministic procedure that sorts all of the $\mathbb{Q}_{\text{cuts}}^{\prime}$-sortable permutations in linear time, since we have such a procedure for the separable permutations \cite{bose}. This is something desirable when considering sorting devices on a restricted class of permutations since the best possible time complexity for a sorting algorithm over all permutations is $\mathcal{O}(n\log{n})$. The popular greedy stack sorting gives such a linear deterministic procedure for stack. The PhD thesis of Luca Ferrari \cite[Section 3.4]{ferrari} shows that such a procedure exists for input-restricted and output-restricted deques, and does not exist for deque. 

Furthermore, popular sorting algorithms, such as Bubblesort, Insertion Sort and Selection Sort correspond to deterministic sorting procedures using certain shuffle queues corresponding to simple shuffling methods.

Sorting by cuts turns out to be an important problem connected to genome rearrangements and an object of extensive study from the algorithms community. For more details, we refer to the introduction of \cite{hartman}. In particular, if we have two permutations representing sequences of genes, we want to find the shortest sequence of operations in a given set that transforms one of the permutations into the other. Assuming that one of the permutations is the identity, the problem is to find the shortest way of sorting a permutation using the fixed set of operations, e.g., cuts and others. The article of Eriksson et al. \cite{shifts} is one work motivated by genome rearrangements that contains results on sorting by cuts which are closest to the bounds we obtain in Theorems \ref{th:nOver2} and \ref{th:lowerB}. They establish bounds for the maximum number of cuts one must apply when sorting a permutation, while we give bounds for the maximum number of iterations of $\mathbb{Q}_{\text{cuts}}^{\prime}$ needed to sort a permutation. The two problems are different, since during an iteration one can apply multiple cuts. Several other articles addressing sorting by cuts together with additional operations, e.g. possible reversions, are listed in \cite{cranston}.

Finally, considering sorting by shuffle queues of type $(ii)$ is reasonable since pop-sorting has been sufficiently considered in the past (see \cite[Chapter 2.1.4]{KitBook}). In addition, we formulate a surprising conjecture involving shuffle queues of type $(ii)$ (see Section \ref{sec:conj}).

\subsection{Summary of our results}

The article is organized as follows.

In Section \ref{sec:deque}, we show that by a deque one can sort the same set of permutations, as with a shuffle queue that can reverse its content. We also show that this is not true for a stack and any given shuffling method.

In Section \ref{sec:cuts}, we study sorting by the device $\mathbb{Q}_{\text{cuts}}^{\prime}$. We show that $S_{n}(\mathbb{Q}_{\text{cuts}}^{\prime})$ is the permutation class $Av_{n}(321,2413,3142)$. A recurrence relation is known for the number of permutations in this class, and thus we get such a relation for $p_{n}(\mathbb{Q}_{\text{cuts}}^{\prime})$. We generalize this result by giving a formula for $p(\mathbb{Q}_{\Sigma}^{\prime})$, for every shuffling method $\Sigma$, such that $\sigma$ is an irreducible permutation for every $\sigma\in \Pi(\Sigma)$, i.e., one for which $\pi([j])\neq [j]$, for any $0<j<n$.

Section \ref{sec:higherCost} investigates permutations having cost greater than one, where $\cost(\pi)$ is the minimal number of times one has to use $\mathbb{Q}_{\text{cuts}}^{\prime}$ in order to sort $\pi$. A natural quantity of interest is $M(n) \coloneqq \max\limits_{\pi\in S_{n}} \cost(\pi)$. We establish bounds from above and below for $M(n)$. As we mentioned, the work of Eriksson et al. \cite{shifts} considers a similar problem and obtain similar bounds for an analogous quantity. The section continues with a conjecture on the limiting behaviour of $M(n)$. We conclude with a proof that $\cost(\pi)=\cost(\pi^{*})$ for every permutation $\pi$, where $\pi^{*}$ is the reverse of the complement of $\pi$.

Section \ref{sec:popDev} is dedicated to pop shuffle queues. First, we prove a statement generalizing the fact that one can sort any given permutation by using $\mathbb{Q}_{\text{cuts}}^{\textsf{pop}}$. In Section \ref{sec:backFront}, we show that the number of permutations of size $n$, sortable by a pop shuffle queue corresponding to any shuffling method with a specific property is enumerated by the Fibonacci numbers $F_{2n-1}$. This fact is an analogue of Theorem \ref{th:irredPerms} for shuffle queues of type $(i)$. Section \ref{sec:conj} discusses a surprising conjecture related to the pop shuffle queues of two popular shuffling methods, namely the \emph{In-shuffle} and the \emph{Monge} shuffling methods. The conjectured fact is that the two methods are \emph{Wilf-pop-equivalent}, that is,  $p_{n}(\mathbb{Q}_{\text{In-sh}}^{\textsf{pop}}) = p_{n}(\mathbb{Q}_{\text{Monge}}^{\textsf{pop}})$ for every $n\geq 1$. We prove that the statement holds if one has to use a single pop operation with each device. Furthermore, we find recursive formulas for the permutations in $S_{n}(\mathbb{Q}_{\text{Monge}}^{\textsf{pop}})$ that end or do not end with $n$, respectively. The same formulas are obtained for $S_{n}(\mathbb{Q}_{\text{In-sh}}^{\textsf{pop}})$ in inequality form and the conjecture holds if and only if these can be replaced by equalities. Using the latter, we have checked that the conjecture holds for $n<20$.

Section \ref{sec:questions} suggests questions for further research.

\section{Shuffle queues equivalent to deque and stack}
\label{sec:deque}
As we explained in Section \ref{sec:motivation}, one motivation to consider shuffle queues is that sorting by deque turns out to be equivalent to sorting by the shuffle queue of a very simple shuffling method that can just reverse its content.
\begin{definition} The shuffling method $\text{rev}$ is defined by the following permutation family.
\begin{equation*}
\label{eq:reverse}
\forall n\geq 2: \Pi_{\text{rev}}^{n} = \{n(n-1)\cdots 21\}.    
\end{equation*}
\end{definition}
For a sequence $w$, the reverse of $w$ will be denoted by $w^{r}$. 
\begin{definition}
The sorting devices $\mathbb{U}$ and $\mathbb{V}$ are \emph{equivalent} if for every $n\geq 1$,
\begin{equation*}
    S_{n}(\mathbb{U}) = S_{n}(\mathbb{V})
\end{equation*}
We denote that by writing $\mathbb{U}\cong \mathbb{V}$.
\end{definition}

\begin{theorem}
\label{th:deque}
$\mathbbm{Deq}\cong \mathbb{Q}_{\text{rev}}$.
\end{theorem}
\begin{proof}
$[$First part: $S_{n}(\mathbbm{Deq})\subseteq S_{n}(\mathbb{Q}_{\text{rev}})]$
Let $\pi\in S_{n}(\mathbbm{Deq})$. Then, there exists an iteration of $\mathbbm{Deq}$ over $\pi$ that sorts it. Take one such iteration $\textsf{itr}$, determined by a sequence of the operations $I,O,\overline{I}$ and $\overline{O}$. Using this sequence, we can easily construct an iteration of $\mathbb{Q}_{\text{rev}}$ that sorts $\pi$, as follows. Instead of the operation $shuffle$ over $\mathbb{Q}_{\text{rev}}$, we will write $\textsf{reverse}$. Replace $I$ by $\textsf{push}$, $O$ by $\textsf{pop}$, $\overline{I}$ by $\textsf{reverse},\textsf{push},\textsf{reverse}$ and $\overline{O}$ by $\textsf{reverse},\textsf{pop},\textsf{reverse}$. This yields a list of operations defining an iteration over $\mathbb{Q}_{\text{rev}}$, which modifies $\pi$ in the exact same way as $\textsf{itr}$ has modified $\pi$ over $\mathbbm{Deq}$.

$[$Second part: $S_{n}(\mathbb{Q}_{\text{rev}}) \subseteq S_{n}(\mathbbm{Deq})]$
If $s$ is a sequence of operations over $\mathbbm{Deq}$, then denote by $\overline{s}$ the \emph{complement} sequence obtained by swapping $I \leftrightarrow \overline{I}$ and $O \leftrightarrow \overline{O}$. Take $\pi\in S_{n}(\mathbb{Q}_{\text{rev}})$ and a sequence of operations $s$ corresponding to an iteration of $\mathbb{Q}_{\text{rev}}$ that sorts $\pi$. The sequence $s$ consists of $\textsf{push}$, $\textsf{pop}$ and $\textsf{reverse}$ operations. Replace every $\textsf{push}$ by an $I$ and every $\textsf{pop}$ by an $O$ to obtain a sequence $s'$. Then, for each $\textsf{reverse}$ operation in $s'$, from left to right, replace the sequence of operations to its left by its complement sequence and then delete that \textsf{reverse} operation. We claim that you will obtain a sequence of operations $s''$ for the device $\mathbbm{Deq}$ that sorts $\pi$. 
For example, suppose that
\[
s= \textsf{push}, \textsf{push}, \textsf{reverse}, \textsf{pop}, \textsf{reverse}, \textsf{push}, \textsf{pop}, \textsf{push}, \textsf{reverse}, \textsf{pop}, \textsf{pop}.
\]
Then,
\[
s'= I, I, \textsf{reverse}, O, \textsf{reverse}, I, O, I, \textsf{reverse}, O, O.
\]
We have three reverse operations in $s'$. If we follow the described procedure, we get:
\[
s' \rightsquigarrow \overline{I}, \overline{I}, O, \textsf{reverse}, I, O, I, \textsf{reverse}, O, O
\]
\[
\rightsquigarrow I, I, \overline{O}, I, O, I, \textsf{reverse}, O, O
\]
\[
\rightsquigarrow \overline{I}, \overline{I}, O, \overline{I}, \overline{O}, \overline{I}, O, O \equalscolon s''.
\]
We will show that the iteration over $\mathbbm{Deq}$ corresponding to $s''$ always sort $\pi$. Assume that $s$ (respectively $s'$) has $r$ reverse operations denoted by $\text{rev}_{i}$ (respectively $\text{rev}'_{i}$), for $i\in [r]$. Furthermore, while transforming $s'$ to $s''$, let the sequence of operations preceding $\text{rev}_{i}$, before replacing it with its complement sequence, be denoted by $s_{(i)}$, for $i\in [r]$. Note that the complement sequence of $s_{(i)}$ is denoted by $s'_{(i)}$, for $i\in [r]$. Our goal is to prove that $s_{(i)}$ transforms $\pi$ in the same way as $s'_{(i)}$, for $i\in [r]$. We will proceed by induction. The sequence $s_{(1)}$ transforms $\pi$ in the same way as $s'_{(1)}$ since $s_{(1)}$ is the complement of $s'_{(1)}$ with a reverse operation added at the end and it is easy to see that if $s$ is a sequence of operations over $\mathbbm{Deq}$ that produces output $\pi'$ on input $\pi$, then $\overline{s}$ produces $(\pi')^{r}$ on input $\pi$. Therefore, if $s'_{(1)}$ produces output $\pi'_{(1)}$ on input $\pi$, then $s_{(1)}$ produces the same output $((\pi'_{(1)})^{r})^{r} = \pi'_{(1)}$ on input $\pi$. Assume that the statement holds for all $i\leq t$ and that $t < r$. By the induction hypothesis, $s'_{(t)}$ transforms the input $\pi$ in the same way as $s_{(t)}$. To obtain $s'_{(t+1)}$ and $s_{(t+1)}$, respectively from $s'_{(t)}$ and $s_{(t)}$, we should first add the same sequence of push and pop operations. Then we take the complement of $s'_{(t)}$ and add a reverse operation to $s_{(t)}$, respectively. We obtain the sequences $s'_{(t+1)}$ and $s_{(t+1)}$ that obviously transform the input $\pi$ in the same way. If $t = r$, then we just add the same sequence of push and pop operations to $s'_{(t)}$ and $s_{(t)}$ to obtain $s$ and $s''$, respectively. Therefore, these two sequences transform $\pi$ in the same way and thus the iteration over $\mathbbm{Deq}$ corresponding to $s''$ also sorts $\pi$.
\end{proof}
Once we know that Theorem \ref{th:deque} holds, a reasonable question to ask is whether there exists a shuffle queue that is equivalent to a stack. Recall that the device stack is denoted by $\mathbbm{St}$.
\begin{theorem}
There is no shuffling method $\Sigma$, such that $\mathbbm{St}\cong \mathbb{Q}_{\Sigma}$.
\end{theorem}
\begin{proof}
Suppose that such a shuffling method $\Sigma$ exists. Then, we must have $S_{n}(\mathbb{Q}_{\Sigma}) = Av_{n}(231)$. Therefore, since $21\in Av_{2}(231)$, we must have $21^{-1}=21\in \Pi_{\Sigma}^{2}$. We also have $231\notin S_{3}(\mathbb{Q}_{\Sigma})$. If $321\in\Pi_{\Sigma}^{3}$, then we will be able to sort $231$ by the following iteration:
\begin{equation*}
\begin{split}
    \begin{pmatrix}
231\\
\varepsilon\\
\varepsilon\\
\end{pmatrix}
\xrightarrow{\textsf{push}}
\begin{pmatrix}
31\\
2\\
\varepsilon\\
\end{pmatrix}
\xrightarrow{\textsf{push}}
\begin{pmatrix}
1\\
23\\
\varepsilon\\
\end{pmatrix}
\xrightarrow[(\text{by }21)]{\substack{\textsf{shuffle}}}
\begin{pmatrix}
1\\
32\\
\varepsilon\\
\end{pmatrix} \\[20pt]
\xrightarrow{\textsf{push}}
\begin{pmatrix}
\varepsilon\\
321\\
\varepsilon\\
\end{pmatrix} 
\xrightarrow[(\text{by }321)]{\substack{\textsf{shuffle}}}
\begin{pmatrix}
\varepsilon\\
123\\
\varepsilon\\
\end{pmatrix}
\xrightarrow{\textsf{pop}}
\begin{pmatrix}
\varepsilon\\
\varepsilon\\
123\\
\end{pmatrix}
\end{split}
\end{equation*}

Thus $321\notin\Pi_{\Sigma}^{3}$. However, we have that $321\in S_{3}(\mathbb{Q}_{\Sigma})$. Consider an input $321$. In order to obtain $123$, a pop operation must not be performed before the first three pushes. Note that after pushing the first two elements, one can either switch them or not, since $21\in \Pi_{\Sigma}^{2}$. Therefore, after pushing the third element $1$, one could either have $231$ or $321$ in the device. Thus we can sort $321$ only if $321^{-1} = 321\in \Pi_{\Sigma}^{3}$ or if $231^{-1} = 312\in \Pi_{\Sigma}^{3}$. However, we saw that $321\notin \Pi_{\Sigma}^{3}$. In addition, $312\notin \Pi_{\Sigma}^{3}$, since otherwise we would be able to sort $231$. This is a contradiction.
\end{proof}
In Section \ref{sec:questions}, we ask a more general question related to shuffle queues equivalent to devices that can sort all the permutations in a given permutation class.
\section{Sorting by cuts}
\label{sec:cuts}
One of the simplest shuffling methods is shuffling by cuts. Its permutation family is given by Equation \eqref{eq:defFcuts}. Some previous works containing results on shuffling using cuts are \cite{diaconisCuts, fulman}. The significance of sorting by $\mathbb{Q}_{\text{cuts}}$ and $\mathbb{Q}_{\text{cuts}}^{\prime}$ is discussed in Section \ref{sec:motivation}. Sorting by $\mathbb{Q}_{\text{cuts}}$ turns out to be trivial since one can sort every given permutation with this shuffle queue. A more general statement is proved at the beginning of Section \ref{sec:popDev}. In this section, we investigate sorting by $\mathbb{Q}_{\text{cuts}}^{\prime}$. Example \ref{ex:cutIter} shows one possible iteration of this device. 

First, we determine $S_{n}(\mathbb{Q}_{\text{cuts}}^{\prime})$, with the help of Lemma \ref{lemma:cutSort}. We will call it the set of the \emph{cut-sortable permutations}. We obtain that this is the set of the separable permutations avoiding the pattern $321$. A permutation $\pi = \pi_{1}\cdots \pi_{n}$ is \emph{separable} if it avoids the patterns $3142$ and $2413$. This important class of permutations arose in the study of pop-stack sorting \cite{avis}. They have a remarkable recursive description and are enumerated by the Schr\"{o}der numbers \cite[Chapter 2.2.5]{KitBook}.

\begin{definition}[Direct sum and skew-sum]
If $\sigma$ and $\tau$ are two permutations of sizes $k$ and $l$, respectively, then their \emph{direct sum} $\sigma\oplus \tau$ and their \emph{skew-sum} $\sigma\ominus \tau$ are defined as follows:
\vspace{-10mm}
\begin{multicols}{2}
  \begin{equation*}
  \Scale[0.9]{
(\sigma\oplus \tau)(i) = \begin{cases}
			\sigma(i), & \text{if }i\leq k\text{,}\\
            k+\tau(i-k), & \text{if $k+1\leq i\leq k+l$.}
		 \end{cases}    
		 }
\end{equation*}\break
\begin{equation*}
\Scale[0.9]{
(\sigma\ominus \tau)(i) = \begin{cases}
			l+\sigma(i), & \text{if $i\leq k$,}\\
            \tau(i-k), & \text{if $k+1\leq i\leq k+l$.}
		 \end{cases}   
		 }
\end{equation*}

\end{multicols}
\end{definition}

\begin{lemma}
\label{lemma:cutSort}
A permutation $\pi$ is in $S_{n}(\mathbb{Q}_{\text{cuts}}^{\prime})$ if and only if it has one of the forms:
\begin{itemize}
    \item[1.] $\pi = id_{r}\oplus\pi'$, for some $1\leq r \leq n$ and $\pi'\in S_{n-r}(\mathbb{Q}_{\text{cuts}}^{\prime})$.
    \item[2.] $\pi = (id_{r_{1}}\ominus id_{r_{2}})\oplus \pi''$, for some $r_{1},r_{2}\geq 1$, where $r\coloneqq r_{1} + r_{2}\leq n$ and $\pi''\in S_{n-r}(\mathbb{Q}_{\text{cuts}}^{\prime})$.
\end{itemize}
\end{lemma}
\begin{proof}
Let $\pi = \pi_{1}\cdots\pi_{n}\in S_{n}(\mathbb{Q}_{\text{cuts}}^{\prime})$. Consider an iteration of $\mathbb{Q}_{\text{cuts}}^{\prime}$ over $\pi$ that sorts it. The sequence of operations for this iteration must contain at least one pop operation. Let the first pop operation be performed after we have pushed $r$ elements in the device ($1\leq r \leq n$), i.e., the elements $\pi_{1},\dots ,\pi_{r}$. The output string after this pop operation must be $id_{r}$. We can have at most one shuffle operation before the first pop operation, and this shuffle must be right before the pop. If we do not have such a shuffle, then the content of the device has not been modified, i.e., $\pi_{1}\cdots\pi_{r} = id_{r}$. Thus, $\pi = id_{r}\oplus \pi'$ and the rest of the iteration sorts $\pi'$. Therefore, $\pi'\in S_{n-r}(\mathbb{Q}_{\text{cuts}}^{\prime})$. If a shuffle has been performed before the first pop, then before this shuffle, the device must contain one of the permutations in the set $(\Pi_{\text{cuts}}^{r})^{-1} = \Pi_{\text{cuts}}^{r}$. Each permutation in $\Pi_{\text{cuts}}^{r}$ can be written as $id_{r_{1}}\ominus id_{r_{2}}$ for some $r_{1},r_{2}\geq 1$, such that $r\coloneqq r_{1}+r_{2} \leq n$. Therefore, $\pi = (id_{r_{1}}\ominus id_{r_{2}})\oplus \pi''$ for some permutation $\pi''\in S_{n-r}(\mathbb{Q}_{\text{cuts}}^{\prime})$  since $\pi''$ is sortable by the rest of the considered iteration. Conversely, one can directly check that any permutation in one of the two listed forms belongs to $S_{n}(\mathbb{Q}_{\text{cuts}}^{\prime})$. 
\end{proof}

 \begin{figure}[ht!]
            \centering
\begin{subfigure}[t]{.5\textwidth}
  \centering
   \begin{tikzpicture}[scale=0.75]
   \filldraw [black] (-2.15,0.35) circle (3pt);
    \filldraw [black] (-1.65,0.85) circle (2.5pt);
    \filldraw [black] (-1.3,1.2) circle (1pt);
    \filldraw [black] (-0.95,1.55) circle (1pt);
    \filldraw [black] (-0.6,1.9) circle (2.5pt);
    
        \draw[dashed] (-3.25,2.5) -- (2.5,2.5);
        \draw[dashed] (0,-0.5) -- (0,5);
        
        \draw (0,2.5) -- (2.5,2.5);
        \draw (0,2.5) -- (0,5);
        \node[below] at (0.85,3.85) {\large{$\pi'$}}; 
    \end{tikzpicture}
  \caption{Cut-sortable permutations that require \\no shuffle before the first pop}
\end{subfigure}%
\begin{subfigure}[t]{.5\textwidth}
  \centering
  \begin{tikzpicture}[scale=0.75]
          \filldraw [black] (-3,1) circle (2.5pt);
    \filldraw [black] (-2.5,1.5) circle (2.5pt);
    \filldraw [black] (-2.25,1.75) circle (1pt);
    \filldraw [black] (-2,2) circle (1pt);
    \filldraw [black] (-1.75,2.25) circle (2.5pt);
    
    \filldraw [black] (-1.5,-0.5) circle (2.5pt);
    \filldraw [black] (-1,0) circle (2.5pt);
    \filldraw [black] (-0.75,0.25) circle (1pt);
    \filldraw [black] (-0.5,0.5) circle (1pt);
    \filldraw [black] (-0.25,0.75) circle (2.5pt);
    
        \draw[dashed] (-3.25,2.5) -- (2.5,2.5);
        \draw[dashed] (0,-0.5) -- (0,5);
        \draw[dashed] (-3.25,0.875) -- (0,0.875);
        \draw[dashed] (-1.625,-0.5) -- (-1.625,2.5);
        
        \draw (0,2.5) -- (2.5,2.5);
        \draw (0,2.5) -- (0,5);
        \node[below] at (0.85,3.85) {\large{$\pi''$}}; 
    
    \end{tikzpicture}
  \caption{Cut-sortable permutations that require \\a shuffle before the first pop}
\end{subfigure}
\caption{}
             \label{fig:lemmaCutSort}
        \end{figure}
An equivalent formulation of Lemma \ref{lemma:cutSort} is that the set $S_{n}(\mathbb{Q}_{\text{cuts}}^{\prime})$ consists of the permutations that can be obtained by direct sums of the trivial permutation $1$ and permutations of the kind $id_{r_{1}}\ominus id_{r_{2}}$. The fact that $S_{n}(\mathbb{Q}_{\text{cuts}}^{\prime})$ is a permutation class follows directly from a simpler version of the observation used to obtain Proposition 1 in \cite{permutingMachines}. With the next theorem, we find this class.
\begin{theorem}
\label{th:cost1}
The permutations sortable by $\mathbb{Q}_{\text{cuts}}^{\prime}$ are the $321$-avoiding separable permutations \cite[A034943]{OEIS}; i.e.,
\begin{equation}
    S_{n}(\mathbb{Q}_{\text{cuts}}^{\prime}) = Av_{n}(321,2413,3142).
\end{equation}
\end{theorem}
\begin{proof}
Let $T \coloneqq \{321,2413,3142\}$. \\
$[$First part: $\pi$ is cut-sortable $\Rightarrow$ $\pi\in Av_{n}(T)]$ We will use induction, Lemma \ref{lemma:cutSort} and the fact that if $\pi = x \oplus y$ for some permutations $x,y$ and $\pi$ has an occurrence of a pattern in $T$, then this occurrence is either in the part of $\pi$ corresponding to $x$ or the part corresponding to $y$. This will be called the \emph{indecomposable property} of $T$. 

The empty permutation belongs to $Av_{0}(T)$. Let $n>0$. Assume, inductively, that any cut-sortable permutation of size $m < n$ belongs to $Av_{m}(T)$. Suppose that $\pi\in S_{n}$ is cut-sortable and $\pi = id_{r}\oplus \pi'$, for some $1\leq r \leq n$ and $\pi'\in S_{n-r}(\mathbb{Q}_{\text{cuts}}^{\prime})$, as in the first form described in Lemma \ref{lemma:cutSort}. Then $\pi'\in Av_{n-r}(T)$ by the inductive hypothesis and $id_{r}$ has no occurrence of a pattern in $T$. Therefore, by the indecomposable property of $T$, we have $\pi\in Av_{n}(T)$.

Now suppose that $\pi$ is in the second form described in the lemma, i.e., that $\pi = (id_{r_{1}}\ominus id_{r_{2}})\oplus \pi''$, for some $r_{1},r_{2}\geq 1$, where $r\coloneqq r_{1} + r_{2}\leq n$ and $\pi''\in S_{n-r}(\mathbb{Q}_{\text{cuts}}^{\prime})$. Then, $\pi''\in Av_{n}(T)$ by the induction hypothesis and one can check easily that $id_{r_{1}}\ominus id_{r_{2}}$ has no occurrence of a pattern in $T$. Because of the indecomposable property of $T$, we must have $\pi\in Av_{n}(T)$.

$[$Second part: $\pi\in Av_{n}(T) \Rightarrow$ $\pi$ is cut-sortable$]$ We will use induction, again. The empty permutation is the only permutation in $Av_{0}(T)$, and it is cut-sortable. Let $n>0$ and $\pi = \pi_{1}\cdots \pi_{n}\in Av_{n}(T)$. Consider the consecutive segment $12\cdots r$ in $\pi$ for the greatest possible value of $r$, where $\pi = \pi_{1}\cdots \pi_{l}12\cdots r\pi_{r+l+1}\cdots \pi_{n}$. If $\pi_{1}\cdots \pi_{l}$ is the empty permutation, then $\pi$ has the first form from Lemma \ref{lemma:cutSort}. If not, then $l\in [1,n-r+1]$ and we will show that $\pi$ has the second form from the lemma. 

First, note that $\pi_{1}\cdots\pi_{l}$ must be increasing to avoid a $321$ pattern in $\pi$. Assume that $\pi_{1}\cdots \pi_{l} \neq (r+1)(r+2)\cdots (r+l)$ and let $u\geq 1$ be minimal, such that $\pi_{u}\neq r+u$. We must have that $u\in [1,l]$, $r+l+1\leq n$, $r+u \in \pi_{r+l+1}\cdots \pi_{n}$ and $\pi_{l}>r+u$. If $u>1$, then $\pi_{1} = r+1$ and $(r+1)\pi_{l}1(r+u)$ would form a $2413$ pattern in $\pi$. Consider $u=1$. Note that $\pi_{r+l+1}\neq r+1$ since $r$ was maximal. In fact, $\pi_{r+l+1}>r+1$. If $\pi_{r+l+1}<\pi_{l}$, then $\pi_{l}\pi_{r+l+1}(r+1)$ would form a $321$ pattern, while if $\pi_{r+l+1}>\pi_{l}$, then $\pi_{l}1\pi_{r+l+1}(r+1)$ would form a $3142$ pattern. Therefore, we must have $\pi_{1}\cdots \pi_{l} = (r+1)(r+2)\cdots (r+l)$ and thus $\pi$ has the second form from Lemma \ref{lemma:cutSort}.
\end{proof}
In \cite{savage}, Martinez and Savage showed that $a_{n} \coloneqq av_{n}(321,2413,3142)$ satisfies
\begin{equation*}
    a_{n} = 3a_{n-1} - 2a_{n-2} + a_{n-3},
\end{equation*}
with initial conditions $a_{1} = 1$, $a_{2} = 2$, $a_{3} = 5$. This is sequence A034943 in the OEIS \cite{OEIS}. The recurrence implies that $a_{n} = \Theta(d^{n})$, where the growth rate $d\approx 2.32$. 

The following theorem gives an alternative way to find the total number of sortable permutations when using cuts. An irreducible permutation $\pi$ is one for which $\pi([j])\neq [j]$ for any $0<j<n$, i.e., the first $j$ elements do not occupy the first $j$ positions. By $IP_{n}$, we denote the set of the irreducible permutations of size $n$. They are enumerated by sequence A003319 in \cite{OEIS}. For example, when $n=3$, the only irreducible permutations are $231$, $312$ and $321$ since they do not have $1$, $12$ or $21$ as a prefix.

\begin{theorem}
\label{th:irredPerms}
If $\Pi_{\Sigma}^{k} \subseteq IP_{k}$ for every $k\geq 2$ and $b_{k} \coloneqq |\Pi_{\Sigma}^{k}|$, then 
\begin{equation}\label{eq:generalize}
    p_{n}(\mathbb{Q}_{\Sigma}^{\prime}) = 1 +  \sum\limits_{\substack{k_{1}+\cdots +k_{l} = n-u \\ k_{i}\geq 2, u\geq 0}}\binom{u+l}{l}\prod\limits_{j=1}^{l}b_{k_{j}}.
\end{equation}
\end{theorem}
\begin{proof}
Recall that a subsequence of consecutive elements $\pi_{a}\cdots\pi_{b}$ is called a segment of $\pi$ and that we denote it by $[a,b]$. When we use $\mathbb{Q}_{\text{cuts}}^{\prime}$, the entire content has to be unloaded after each shuffle and the segments of the input that were not shuffled are kept the same in the output. Thus the output after an iteration of $\mathbb{Q}_{\Sigma}^{\prime}$ is uniquely determined by the segments of the input that were shuffled and the corresponding permutations chosen for each of the shuffle operations. For instance, the output $id_{6}$ of the iteration of $\mathbb{Q}_{\text{cuts}}^{\prime}$ shown in Example \ref{ex:cutIter} is determined by the sequence of segments $([1,2],[4,6])$ of the input $213645$ that were shuffled and the sequence of permutations $(21,231)$ that were applied on the given segments. 

Denote the set of the possible pairs of sequences of segments and permutations, for an input of size $n$ and a shuffling method $\Sigma$, by $SSP_{n}^{\Sigma}$. For every $n\geq 2$ and every element $(s,q)\in SSP_{n}^{\Sigma}$, the segments in $s$ are in lexicographical order and do not overlap with each other since we shuffle these segments from left to right. 
We will first show that $|SSP_{n}^{\Sigma}|$ is equal to the expression in the right-hand side of \eqref{eq:generalize}. Then, we will give a bijection between the sets $S_{n}(\mathbb{Q}_{\Sigma}^{\prime})$ and $SSP_{n}^{\Sigma}$.

$[$Finding $|SSP_{n}^{\Sigma}|]$ Assume that $x=(s,q)\in SSP_{n}^{\Sigma}$ and that $s$ consists of $l$ shuffled segments. Only one such $x$ exists, if $l=0$. Let $l\geq 1$. Denote the sizes of the $l$ shuffled segments by $k_{1},k_{2},\dots ,k_{l}$, where $k_{j}\geq 2$ for every $j\in [1,l]$ and let their sum be $n-u$ for some $u\geq 0$. For instance, if $n=8$, $l=2$ and $s = ([2,3],[5,7])$, then $k_{1} = 2$, $k_{2} = 3$ and $n-u = 5$. In general, if the numbers $l,n-u, k_{1},\dots ,k_{l}$ are given, then in order to determine the sequence of segments $s$, one should distribute the $u$ remaining elements in the set of $l+1$ spaces - one before each of the $l$ segments and the one after all of the segments. For every such choice, we obtain a different sequence of segments $s$. The number of these choices is the number of ways to distribute $u$ indistinguishable balls into $l+1$ boxes that is $\binom{u+(l+1)-1}{(l+1)-1} = \binom{u+l}{l}$. Then, if $q = (q_{1},\ldots ,q_{l})$, the permutation $q_{j}$ can be any of the $b_{k_{j}}$ permutations in $\Pi_{\Sigma}^{k_{j}}$ for every $j\in [1,l]$. Thus $q$ can be determined in $\prod\limits_{j=1}^{l}b_{k_{j}}$ ways. In total, we obtain the right-hand side of \eqref{eq:generalize}.

$[S_{n}(\mathbb{Q}_{\Sigma}^{\prime})$ $\rightarrow$ $SSP_{n}^{\Sigma}]$ 
Let $\pi\in S_{n}(\mathbb{Q}_{\Sigma}^{\prime})$. Then, there exists at least one iteration that sorts $\pi$. Assume that $\pi$ can be sorted by two different iterations $it_{1}$ and $it_{2}$, corresponding to $x_{1},x_{2}\in SSP_{n}^{\Sigma}$, where $x_{1}=(s_{1},q_{1})$, $x_{2}=(s_{2},q_{2})$ and $x_{1}\neq x_{2}$. Assume that $s_{1}=s_{2}$. Then, $q_{1}\neq q_{2}$. However, we can easily see that this is not possible. Let $[r,r+k]$ be an arbitrary segment in $s_{1}$, and respectively in $s_{2}$. If $\sigma_{1}$ and $\sigma_{2}$ are the two permutations in $q_{1}$ and $q_{2}$, respectively, that have to be applied on this segment, then we must have $\sigma_{1}=\sigma_{2} = (\pi_{r}\cdots\pi_{r+k})^{-1}$. Thus $q_{1}= q_{2}$. Therefore, we must have $s_{1}\neq s_{2}$. 

Let $[r_{1},r_{1}+k_{1}]$ be the last segment in $s_{1}$ and let $[r_{2},r_{2}+k_{2}]$ be the last segment in $s_{2}$. Assume also that $\sigma_{1}$ and $\sigma_{2}$ are the last permutations in $q_{1}$ and $q_{2}$, respectively. We saw that if $[r_{1},r_{1}+k_{1}] = [r_{2},r_{2}+k_{2}]$, then we must have $\sigma_{1}=\sigma_{2}$. 
However, $s_{1}\neq s_{2}$. Therefore, without loss of generality, assume that $[r_{1},r_{1}+k_{1}] \neq [r_{2},r_{2}+k_{2}]$ and that $r_{1}\leq r_{2}$. If $r_{1} = r_{2}$, then assume for concreteness that $k_{1}<k_{2}$ (see Figure \ref{fig:2segmentsEq}). 

    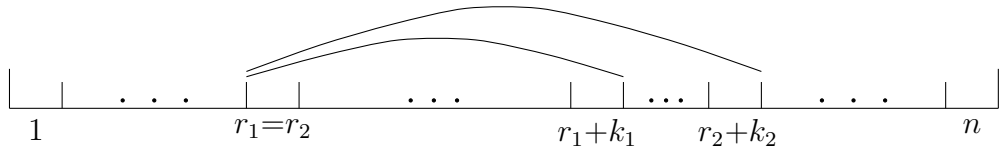
\begin{figure}[h!]
            \centering
\begin{tikzpicture}[scale=0.7]
\draw (0,0) -- (18.78,0);
        \draw (0,0) -- (0,0.75);
        \node[below] at (0.5,0) {$1$}; 
        \draw (1,0) -- (1,0.5);
        
        \filldraw[fill=black] (2.15,0.15) circle [radius=0.03];
        \filldraw[fill=black] (2.75,0.15) circle [radius=0.03];
        \filldraw[fill=black] (3.35,0.15) circle [radius=0.03];
        
        \draw (4.5,0) -- (4.5,0.5);
        \node[below] at (5,0) {$r_{1}{=} r_{2}$}; 
        \draw (5.5,0) -- (5.5,0.5);
        
        \draw [black] plot [smooth] coordinates {(4.5,0.6) (5.95,1) (7.4,1.3) (8.85,1.3) (10.3,1)  (11.66,0.6)};
        
        \draw [black] plot [smooth] coordinates {(4.5,0.7) (5.89,1.2) (7.28,1.6) (8.67,1.9) (10.06,1.9) 
        (11.45,1.6) 
        (12.84,1.2) 
        (14.28,0.7)};
        
        \filldraw[fill=black] (7.62,0.15) circle [radius=0.03];
        \filldraw[fill=black] (8.06,0.15) circle [radius=0.03];
        \filldraw[fill=black] (8.5,0.15) circle
        [radius=0.03];
        
        \draw (10.66,0) -- (10.66,0.5);
        \node[below] at (11.16,0) {$r_{1}{+}k_{1}$};
        \draw (11.66,0) -- (11.66,0.5);
        
        \filldraw[fill=black] (12.19,0.15) circle [radius=0.03];
        \filldraw[fill=black]
        (12.47,0.15) circle [radius=0.03];
        \filldraw[fill=black] (12.75,0.15) circle [radius=0.03];
        
        \draw (13.28,0) -- (13.28,0.5);
        \node[below] at (13.83,0) {$r_{2}{+}k_{2}$};
        \draw (14.28,0) -- (14.28,0.5);
        
        \filldraw[fill=black] (15.43,0.15) circle [radius=0.03];
        \filldraw[fill=black]
        (16.03,0.15) circle [radius=0.03];
        \filldraw[fill=black] (16.63,0.15) circle [radius=0.03];
        
        \draw (17.78,0) -- (17.78,0.5);
            \node[below] at (18.28,0) {$n$};
        \draw (18.78,0) -- (18.78,0.75);
    \end{tikzpicture}
    \caption{The case $r_{1}=r_{2}$}
    \label{fig:2segmentsEq}
\end{figure}

Iteration $it_{1}$ permutes the elements of the segment $[r_{1},r_{1}+k_{1}]$ in $\pi$. Hence $it_{2}$ does the same. This implies that $\sigma_{2}$ fixes $[k_{1}+1]$ and thus $\sigma_{2}\notin IP_{k_{2}+1}$, which is a contradiction. If $r_{1} < r_{2}$, then it suffices to look at the following two cases (see Figure \ref{fig:2segments1} and Figure \ref{fig:2segments2}):
\begin{enumerate}
    \item $r_{2}\leq r_{1}{+}k_{1}$.
    Then, $\sigma_{1}$ fixes $[r_{2}-r_{1}]$. Indeed, suppose that $\sigma_{1}(u)=v$, where $u\in [r_{2}-r_{1}]$ and $v>r_{2}-r_{1}$. This means that $it_{1}$ moves $\pi_{r_{1}+u-1}$ to position $r_{1}+v-1\geq r_{1}+(r_{2}-r_{1}+1)-1 = r_{2}$. However, $it_{2}$ moves $\pi_{r_{1}+u-1}$ to a position smaller than $r_{2}$. Therefore, $\sigma_{1}$ fixes $[r_{2}-r_{1}]$ and $\sigma_{1}$ is not irreducible, which is a contradiction.
    \begin{figure}[h!]
            \centering
\begin{tikzpicture}[scale=0.7]
\draw (0,0) -- (18.78,0);
        \draw (0,0) -- (0,0.75);
        \node[below] at (0.5,0) {$1$}; 
        \draw (1,0) -- (1,0.5);
        
        \filldraw[fill=black] (2.15,0.15) circle [radius=0.03];
        \filldraw[fill=black] (2.75,0.15) circle [radius=0.03];
        \filldraw[fill=black] (3.35,0.15) circle [radius=0.03];
        
        \draw (4.5,0) -- (4.5,0.5);
        \node[below] at (5,0) {$r_{1}$}; 
        \draw (5.5,0) -- (5.5,0.5);
        
        \draw [black] plot [smooth] coordinates {(4.5,0.6) (5.95,1) (7.4,1.4) (8.85,1.4) (10.3,1)  (11.66,0.6)};
        
        \filldraw[fill=black] (6.03,0.15) circle [radius=0.03];
        \filldraw[fill=black] (6.31,0.15) circle [radius=0.03];
        \filldraw[fill=black] (6.59,0.15) circle
        [radius=0.03];
        
        \draw (7.12,0) -- (7.12,0.5);
        \node[below] at (7.57,0) {$r_{2}$}; 
        \draw (8.12,0) -- (8.12,0.5);
        
        \draw [black] plot [smooth] coordinates {(7.12,0.6) (8.57,1) (10.02,1.4) (11.47,1.4) (12.92,1)  (14.28,0.6)};
        
        \filldraw[fill=black] (8.96,0.15) circle [radius=0.03];
        \filldraw[fill=black] (9.4,0.15) circle [radius=0.03];
        \filldraw[fill=black] (9.84,0.15) circle
        [radius=0.03];
        
        \draw (10.66,0) -- (10.66,0.5);
        \node[below] at (11.16,0) {$r_{1}{+}k_{1}$};
        \draw (11.66,0) -- (11.66,0.5);
        
        \filldraw[fill=black] (12.19,0.15) circle [radius=0.03];
        \filldraw[fill=black]
        (12.47,0.15) circle [radius=0.03];
        \filldraw[fill=black] (12.75,0.15) circle [radius=0.03];
        
        \draw (13.28,0) -- (13.28,0.5);
        \node[below] at (13.83,0) {$r_{2}{+}k_{2}$};
        \draw (14.28,0) -- (14.28,0.5);
        
        \filldraw[fill=black] (15.43,0.15) circle [radius=0.03];
        \filldraw[fill=black]
        (16.03,0.15) circle [radius=0.03];
        \filldraw[fill=black] (16.63,0.15) circle [radius=0.03];
        
        \draw (17.78,0) -- (17.78,0.5);
            \node[below] at (18.28,0) {$n$};
        \draw (18.78,0) -- (18.78,0.75);
    \end{tikzpicture}
    \caption{The case $r_{1}<r_{2}$ and $r_{2}\leq r_{1}+k_{1}$}
    \label{fig:2segments1}
\end{figure}
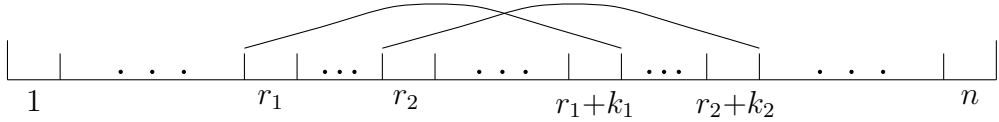

\item $r_{2}>r_{1}{+}k_{1}$. 
    Since $it_{1}$ sorts $\pi$, we must have $\sigma_{2}=id_{k_{2}+1}$, which is not possible.
    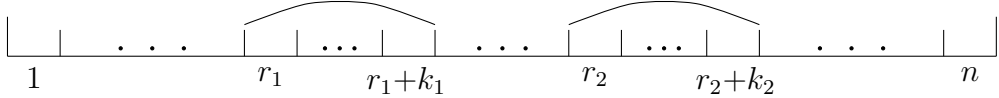
\begin{figure}[H]
            \centering
\begin{tikzpicture}[scale=0.7]
\draw (0,0) -- (18.78,0);
        \draw (0,0) -- (0,0.75);
        \node[below] at (0.5,0) {$1$}; 
        \draw (1,0) -- (1,0.5);
        
        \filldraw[fill=black] (2.15,0.15) circle [radius=0.03];
        \filldraw[fill=black] (2.75,0.15) circle [radius=0.03];
        \filldraw[fill=black] (3.35,0.15) circle [radius=0.03];
        
        \draw (4.5,0) -- (4.5,0.5);
        \node[below] at (5,0) {$r_{1}$}; 
        \draw (5.5,0) -- (5.5,0.5);
        
        \draw [black] plot [smooth] coordinates {(4.5,0.6) (5.71,1) (6.92,1) (8.12,0.6)};
        
        \filldraw[fill=black] (6.03,0.15) circle [radius=0.03];
        \filldraw[fill=black] (6.31,0.15) circle [radius=0.03];
        \filldraw[fill=black] (6.59,0.15) circle
        [radius=0.03];
        
        \draw (7.12,0) -- (7.12,0.5);
        \node[below] at (7.57,0) {$r_{1}{+}k_{1}$}; 
        \draw (8.12,0) -- (8.12,0.5);
        
        \filldraw[fill=black] (8.96,0.15) circle [radius=0.03];
        \filldraw[fill=black] (9.4,0.15) circle [radius=0.03];
        \filldraw[fill=black] (9.84,0.15) circle
        [radius=0.03];
        
        \draw (10.66,0) -- (10.66,0.5);
        \node[below] at (11.16,0) {$r_{2}$};
        \draw (11.66,0) -- (11.66,0.5);
        
        \draw [black] plot [smooth] coordinates {(10.66,0.6) (11.87,1) (13.08,1) (14.28,0.6)};
        
        \filldraw[fill=black] (12.19,0.15) circle [radius=0.03];
        \filldraw[fill=black]
        (12.47,0.15) circle [radius=0.03];
        \filldraw[fill=black] (12.75,0.15) circle [radius=0.03];
        
        \draw (13.28,0) -- (13.28,0.5);
        \node[below] at (13.83,0) {$r_{2}{+}k_{2}$};
        \draw (14.28,0) -- (14.28,0.5);
        
        \filldraw[fill=black] (15.43,0.15) circle [radius=0.03];
        \filldraw[fill=black]
        (16.03,0.15) circle [radius=0.03];
        \filldraw[fill=black] (16.63,0.15) circle [radius=0.03];
        
        \draw (17.78,0) -- (17.78,0.5);
            \node[below] at (18.28,0) {$n$};
        \draw (18.78,0) -- (18.78,0.75);
    \end{tikzpicture}
    \caption{The case $r_{1}<r_{2}$ and $r_{2}> r_{1}+k_{1}$}
    \label{fig:2segments2}
\end{figure}
\end{enumerate}
We see that it is not possible to sort $\pi$ by two different iterations corresponding to two different elements of $SSP_{n}^{\Sigma}$. Therefore, for every $\pi\in S_{n}(\mathbb{Q}_{\Sigma}^{\prime})$ there exists a unique $x\in SSP_{n}^{\Sigma}$ corresponding to an iteration that sorts $\pi$.

$[SSP_{n}^{\Sigma}$ $\rightarrow$ $S_{n}(\mathbb{Q}_{\Sigma}^{\prime})]$
It remains to show that every $x\in SSP_{n}^{\Sigma}$ corresponds to a set of iterations of $\mathbb{Q}_{\Sigma}^{\prime}$ sorting exactly one permutation $\pi$. Let $x = (s,q)$, where $s = (s_{1},\ldots ,s_{l})$ and $q = (\sigma_{1}, \ldots ,\sigma_{l})$. Take $id_{n}$, and go backwards by applying consecutively $\sigma_{j}^{-1}$ to the segment $s_{j}$, for $j= l,l-1, \ldots ,1$. We will obtain a unique permutation $\pi$ that is sortable by any iteration $\textsf{itr}$ corresponding to $x$.
\end{proof}

Note that when $\Sigma = cuts$, we have $\Pi_{\text{cuts}}^{k} \subseteq IP_{k}$ and $b_{k} = |\Pi_{\text{cuts}}^{k}| = k-1$, for every $k\geq 2$. Thus, one can apply formula \eqref{eq:generalize} in order to compute $p_{n}(\mathbb{Q}_{\text{cuts}}^{\prime})$.

\section{Permutations of higher cost}
\label{sec:higherCost}
Obviously, not all $\pi\in S_{n}$ are sortable by $\mathbb{Q}_{\text{cuts}}^{\prime}$. However, one can use a device several times in a row by using the output after one iteration as an input to the next iteration. This is the so-called \emph{sorting in series}. Many articles investigate this idea for stack-sorting (see \cite[Section 8.2.2]{Bona}). Denote the set of permutations that one can obtain after $k$ iterations of $\mathbb{Q}_{\text{cuts}}^{\prime}$ over a permutation $\pi\in S_{n}$ by $(\mathbb{Q}_{\text{cuts}}^{\prime})^{k}(\pi)$. 
\begin{definition}[cost of permutation]
The \emph{cost} of $\pi$ is the minimum number of iterations needed to sort $\pi$ using the device $\mathbb{Q}_{\text{cuts}}^{\prime}$, i.e.,

\[\cost(\pi) \coloneqq \min \{m \mid id_{n}\in (\mathbb{Q}_{\text{cuts}}^{\prime})^{m}(\pi)\}.
\]
\end{definition}

It is not difficult to obtain an upper bound for $\cost(\pi)$. Indeed, one can move a single element to its correct position using only one iteration. In particular, if the input permutation is $\pi = 12\cdots (i-1)\pi_{i}\cdots\pi_{j-1}i\pi_{j+1}\cdots\pi_{n}$, then one can perform an iteration consisting of only one cut right before $i$, after getting the subsequence $\pi_{i}\cdots\pi_{j-1}i$ into the device. Such an iteration will move $i$ at its correct position. Consecutive movements of $i,i{+}1,\dots ,n$ to their correct positions will sort the permutation. Therefore, $\cost(\pi)\leq n$. This upper bound is improved significantly with the theorem given below.

\begin{theorem}
\label{th:nOver2}
$\cost(\pi)\leq \lceil \frac{n}{2}\rceil$, for every $\pi\in S_{n}$, where $n\geq 1$.
\end{theorem}
\begin{proof}
A computer simulation shows that the statement is true for $1\leq n \leq 10$. Let $n\geq 11$ and let us assume, inductively, that the statement holds for all $n'<n$. The main observation that will be used is that if we have $k+1$ consecutive numbers in $[n]$, forming a segment in $\pi$, then we can treat them as a single element and apply the induction hypothesis for $n-k$. Two more observations will be substantially used that describe cases when we can modify $\pi$ with one iteration over $\mathbb{Q}_{\text{cuts}}^{\prime}$ and then use the main observation above. A third observation for the case when $\pi$ is a direct sum of other permutations will be also needed. These three observations are listed below with a brief justification for each of them:
\begin{enumerate}
    \item If $a\in [3,n]$ and the numbers $a-1$, $a-2$ occur before $a$ in $\pi$, then there exists $\pi'\in\mathbb{Q}_{\text{cuts}}^{\prime}(\pi)$, such that $\pi' = \dots (a-2)(a-1)a \dots$.
    
    Proof: Assume that $\pi = \pi_{1}\cdots(a-1)\cdots (a-2)\pi_{h}\cdots a\cdots\pi_{n}$ for some $h>2$. If $(a-2)$ is before $(a-1)$, then we can proceed in a similar way. We can perform the following two cuts with one iteration of $\pi$ over $\mathbb{Q}_{\text{cuts}}^{\prime}$: Cut the segment $(a-1)\cdots (a-2)$ after $(a-1)$ and the segment $\pi_{h}\cdots a$ before $a$. A permutation $\pi'$ with the desired property is obtained. If $a=\pi_{h}$, then we will not need the second cut.
    
    \item Assume that $a,a+1,b,b+1\in [n]$ are four different numbers. If $a$ and $a+1$ occur before $b$ and $b+1$ in $\pi$, then there exists $\pi'\in\mathbb{Q}_{\text{cuts}}^{\prime}(\pi)$, such that $\pi' = \ldots a(a+1)\cdots b(b+1)\ldots$
    
    Proof: We can simply move $a+1$ immediately after $a$ with a single cut and $b+1$ immediately after $b$ with another cut. Since $a$ and $a+1$ occur before $b$ and $b+1$, we can perform the two cuts within one iteration.
    
    \item Assume that $\pi = \sigma_{1}\oplus\sigma_{2}\oplus\cdots\oplus\sigma_{k}$.  Then, 
    \[
    \cost(\pi)\leq \max(\cost(\sigma_{1}),\cost(\sigma_{2}),\dots ,\cost(\sigma_{k})).
    \]
    
    Proof: Within one iteration, one may independently transform each of the parts of $\pi$ corresponding to $\sigma_{1},\sigma_{2}, \dots , \sigma_{k}$. Thus, if $m = \max(\cost(\sigma_{1}), \dots ,\cost(\sigma_{k}))$, then $\pi$ can be sorted with $m$ iterations.
\end{enumerate}
We continue with the proof. Let $xy$ denote the last two elements of $\pi$. Observation $(1)$ implies that unless $y = 1$ or $y = 2$, we will be able to transform $\pi$ to a permutation $\pi'$ containing the segment $(y-2)(y-1)y$ with just one iteration. Looking at this segment as a single element and applying the induction hypothesis for $n-2$ would give us $\cost(\pi)\leq 1 + \lceil \frac{n-2}{2}\rceil = \lceil \frac{n}{2}\rceil$, which is what we want. Assume that $y=1$. Obviously, if $x\neq 2,3$, we will be able to apply observation $(1)$ again with $a=x$ and obtain the bound via the same calculation. If $x=2$, then since $n\geq 4$, the numbers $3$ and $4$ will precede $1$ and $2$ in $\pi$ which allows us to use observation $(2)$ and obtain the permutation $\pi'$ described there, with one iteration of $\pi$ over $\mathbb{Q}_{\text{cuts}}^{\prime}$. Treating both $1,2$ and $3,4$ as a single element and applying the induction hypothesis gives us the same calculation and implies the desired result, again. 

Therefore, $xy = 31$ is the only case that remains to be considered if $y=1$. If $y=2$, then if $x=1$, we can move the last two elements $xy=12$ at the beginning of $\pi$ with an iteration consisting of a single cut to obtain a permutation $\pi'=12\oplus\pi''$. Applying observation $(3)$ to $\pi'$ and the induction hypotheses for $\pi''$ gives $\cost(\pi)\leq 1 + \cost(\pi'') \leq 1 + \lceil \frac{n-2}{2}\rceil \leq \lceil \frac{n}{2}\rceil$. Therefore, we may assume that if $y=2$, then $x\neq 1$. If $x\neq 3,4$, then we can obtain the result using observation $(1)$, as before. If $x=3$, then since $n\geq 5$, we will be able to apply observation $(2)$ for $2,3$ and $4,5$. Therefore, $xy = 42$ is the only case that remains to be considered if $y=2$. 

We saw that it suffices to look at those permutations $\pi$ having last two elements, $xy=31$ or $xy=42$. Following the same reasoning, we can easily obtain that it suffices to only look at permutations $\pi$ beginning either with $n(n-2)$ or $(n-1)(n-3)$. The only difference is that an observation analogous to $(1)$ shall be used dealing with the cases when $a-2$ precedes both $a-1$ and $a$. Hence, we have four cases, in total. We will show how we can complete the proof in only one of them, namely when $\pi$ begins with $n(n-2)$ and finishes with $42$. The proofs in the other 3 cases can be completed following the same reasoning.

Let $\pi = n(n-2)\cdots 42$. Then, we can assume that $n-3$ and $n-1$ occur after both $1$ and $3$, because otherwise we will be able to apply observation $(2)$ for certain pairs of elements. For concreteness, let us take $3$ to be before $1$ and $n-1$ to be before $n-3$. The following argument works regardless of this order. We can assume that $\pi = n(n-2)\cdots 3\cdots 1\cdots (n-1)\cdots (n-3)\cdots42$. Below, we show a particular way to transform $\pi$ by four iterations. We give the output at the end of each iteration. The reader may try to find the exact cuts applied in these iterations.
\begin{figure}[h!]
            \centering
\begin{tikzpicture}[scale=0.7]
        \node[below] (A) at (0,10) {$n(n-2)\cdots 3\cdots 1\cdots (n-1)\cdots (n-3)\cdots42$}; 
        \node[below] (B) at (0,8.5) {$n(n-2)\cdots (n-1)\cdots (n-3)\cdots 3\cdots 142$}; 
        \node[below] (C) at (0,7) {$\cdots (n-1)n(n-2)(n-3)\cdots 3142$};
        \node[below] (D) at (0,5.5) {$3142\cdots (n-1)n(n-2)(n-3)\cdots$};
        \node[below] (E) at (0,4) {$3142\cdots (n-1)n(n-2)(n-3)$};
        
        \draw [->] (A) edge (B) (B) edge (C) (C) edge (D) (D) edge (E);
    \end{tikzpicture}
\end{figure}

Therefore, by four iterations $\pi$ can be transformed to $\pi' = w_{1}\pi''w_{2}$, where $|\pi''| = n-8$, $|w_{1}| = 4$, $|w_{2}| = 4$ and $\pi' = \red(w_{1})\oplus \red(\pi'')\oplus \red(w_{2})$. Recall that $n\geq 11$ and thus $|\pi''| = n-8\geq 3$, which means that $\cost(\pi'')\geq 2$. In addition, $\cost(\sigma)=2$, for any $\sigma\in S_{4}$. Therefore, observation $(3)$ applied over $\pi'$ and the induction hypothesis for $\pi''$ gives us $\cost(\pi')\leq \cost(\pi'') \leq \lceil \frac{n-8}{2}\rceil$, which implies $\cost(\pi)\leq 4 + \cost(\pi') \leq 4 + \lceil \frac{n-8}{2}\rceil = \lceil \frac{n}{2}\rceil$.

\end{proof}

Theorem \ref{th:nOver2} gives a tight upper bound for the $\cost$ function since there exist permutations of size $n$ and cost $\lceil \frac{n}{2}\rceil$. For instance, we computed that $\cost(83527461) = 4$. The best absolute lower bound is obviously $0$ since $\cost(id_{n})=0$, for every $n$. Let $M(n) \coloneqq \max\limits_{\pi\in S_{n}} \cost(\pi)$ be the maximal cost of a permutation of size $n$. Theorem \ref{th:nOver2} gives us that $M(n)\leq \lceil \frac{n}{2}\rceil$. Next, we give a lower bound for $M(n)$ by Theorem \ref{th:lowerB}. We begin by showing that cost is monotonically increasing with respect to pattern containment. We will write $C_{n}(q)\coloneqq S_{n}\setminus Av_{n}(q)$ for the permutations of $[n]$ that contain the pattern $q$. A main fact that will be used is that sorting by cuts has the property defined below.

\begin{definition}[hereditary property]
A shuffling method $\Sigma$ has the \emph{hereditary} property if the following holds: Suppose that a sequence $\sigma$ can be transformed to a sequence $\sigma'$ by a permutation in $\Pi(\Sigma)$. If $\tau$ is a subsequence of $\sigma$ and its symbols transform to the subsequence $\tau'$ of $\sigma'$, then there exists a permutation in $\Pi(\Sigma)$ transforming $\tau$ to $\tau'$.
\end{definition}
This property is defined in \cite{permutingMachines}, as a property of the so-called \say{permuting machines}. Here we will use that shuffling by cuts has this property.
\begin{lemma}
\label{lemma:pattern}
If $\pi\in C_{n}(q)$, then $\cost(\pi)\geq \cost(q)$.
\end{lemma}
\begin{proof}
Let us fix an occurrence $\textsf{oc}$ of $q$ in $\pi$. Assume that we have a sequence of iterations sorting $\pi$ and let $\textsf{itr}$ be one of these iterations. Every cut $c$ in $\textsf{itr}$ is transforming a certain sequence of elements $\sigma$ to a sequence $\sigma'$. If $\tau$ is the subsequence of $\sigma$, including all of the elements of $\textsf{oc}$, that is transformed to a sequence $\tau'$, then by the hereditary property of sorting by cuts, there exists a cut $c'$ which transforms $\tau$ to $\tau'$. Therefore, for every sequence of iterations of $\mathbb{Q}_{\text{cuts}}^{\prime}$ that sorts $\pi$, one can get a corresponding sequence of iterations that sorts its subsequence $\textsf{oc}$ by substituting each cut $c$ in an iteration from the initial sequence with the corresponding cut $c'$. The total number of iterations may drop since some of the iterations in the initial sequence sorting $\pi$ may not affect the elements of $\textsf{oc}$. If we consider an optimal sequence of $\cost(\pi)$ iterations sorting $\pi$, then the described correspondence gives a sequence of at most $\cost(\pi)$ iterations of $\mathbb{Q}_{\text{cuts}}^{\prime}$ sorting $q$. Thus $\cost(q)\leq \cost(\pi)$.
\end{proof}

The last lemma shows that the cost function is monotone over the partially ordered set of permutations ordered by pattern containment. 

Recall that $id_{n}^{r}$ is the reverse identity: $id_{n}^{r} = n(n-1)\cdots 1$.

\begin{lemma}
\label{lemma:decrSeq}
If $\pi'\in\mathbb{Q}_{\text{cuts}}^{\prime}(id_{n}^{r})$, then $\pi'\in C_{n}(id_{\lceil \frac{n}{2}\rceil}^{r} )$, i.e., $\pi'$ contains a decreasing subsequence of size $\lceil \frac{n}{2}\rceil$.
\end{lemma}
\begin{proof}
We will proceed by induction. The lemma holds for $n=1$. Consider the first cut $c$ in an arbitrary iteration of $\mathbb{Q}_{\text{cuts}}^{\prime}$ over $id_{n}^{r}$. Denote the output permutation after this iteration by $\pi' = \pi_{1}'\cdots\pi_{n}'$. If $n$ does not participate in $c$, then $\pi_{1}' = n$, because we have just pushed and popped $\pi_{1}=n$. Then we can look at the considered iteration as one over $id_{n-1}^{r}$ with the element $n$ appended in front of the output. The element $n$ is in front of any decreasing subsequence in $\pi_{2}'\cdots\pi_{n}'$. Therefore, we may apply the induction hypothesis to get that $\pi'$ must have a decreasing subsequence of size $1 + \lceil \frac{n-1}{2}\rceil \geq \lceil \frac{n}{2}\rceil$. 

If $n$ participates in $c$, then note that a cut can be performed only if we have more than one number in the device. Let the cut $c$ be performed after we have exactly $k\geq 2$ numbers in the device $\mathbb{Q}_{\text{cuts}}^{\prime}$: $n(n-1)\cdots (n-k+1)$. After the cut $c$, these $k$ numbers will be divided into two decreasing sequences. In other words, $\pi_{1}'\cdots\pi_{k}'$ will be comprised of two segments that are decreasing sequences. At least one of these sequences must be of size at least $\lceil \frac{k}{2}\rceil$ and therefore $\pi_{1}'\cdots\pi_{k}'$ contains a decreasing sequence of such size. The rest of the iteration can be looked at as an iteration over $id_{n-k}^{r}$. Thus we can apply the inductive hypothesis to see that $\pi_{k+1}'\cdots\pi_{n}'$ must contain a decreasing sequence of size $\lceil \frac{n-k}{2}\rceil$. In addition, $\pi_{1}'\cdots\pi_{k}' >\pi_{k+1}'\cdots\pi_{n}'$, so $\pi'$ must contain a decreasing sequence of size $\lceil \frac{k}{2}\rceil + \lceil \frac{n-k}{2}\rceil \geq \lceil \frac{n}{2}\rceil$.
\end{proof}
Now, we are ready to establish the lower bound for $M(n)$, i.e., the maximal cost of a permutation of size $n$.
\begin{theorem}
\label{th:lowerB}
$M(n)\geq \lceil \log_{2}{n} \rceil$, for each $n\geq 2$.
\end{theorem}
\begin{proof}
We will prove that $\cost(id_{n}^{r})\geq \lceil \log_{2}{n} \rceil$ for each $n\geq 2$, using induction. Obviously, $id_{2}^{r} = 21$ cannot be sorted with less than one iteration through $\mathbb{Q}_{\text{cuts}}^{\prime}$. By Lemma \ref{lemma:decrSeq}, after one iteration over $id_{n}^{r}$, we will always get a permutation $\pi'\in C_{n}(id_{\lceil \frac{n}{2}\rceil}^{r})$. By Lemma \ref{lemma:pattern} and the induction hypothesis, 
\[
\cost(\pi')\geq \cost(id_{\lceil \frac{n}{2}\rceil}^{r})\geq \left\lceil \log_{2}{\left\lceil \frac{n}{2}\right\rceil} \right\rceil \geq \left\lceil \log_{2}{\frac{n}{2}} \right\rceil = \lceil \log_{2}{n} \rceil -1.
\]
But, $\pi'\in\mathbb{Q}_{\text{cuts}}^{\prime}(id_{n}^{r})$, so $\cost(id_{n}^{r}) = 1 + \cost(\pi')\geq 1 + (\lceil \log_{2}{n} \rceil -1) = \lceil \log_{2}{n} \rceil$.
\end{proof}

There exist values of $n$ for which $M(n)>\lceil \log_{2}{n} \rceil$. We believe that the set of these values is bounded, but we were not able to prove that.
\begin{question}
Is it true that $M(n) \xrightarrow[n\to \infty]{} \lceil \log_{2}{n} \rceil$?
\end{question}

A positive answer to this question would imply that every permutation of size $n$ can be sorted using $\mathcal{O}(n\log n)$ operations by using cuts since one iteration uses $\mathcal{O}(n)$ operations.

We finish this section by showing that the permutations in $S_{n}$ can be paired up in terms of cost, when using $\mathbb{Q}_{\text{cuts}}^{\prime}$. For a permutation $\pi = \pi_{1}\cdots\pi_{n}$, let $\overline{\pi}$ denote the complement permutation, defined by $\overline{\pi}_{i} = n+1 - \pi_{i}$. Recall that $\pi^{r}$ denotes the reverse of $\pi$, i.e., $(\pi^{r})_{i} = \pi_{n+1-i}$. Set $\pi^{*} = \overline{\pi^{r}} = (\overline{\pi})^{r}$. Observe also that $\Pi_{\text{cuts}}^{n}$ is closed under the $^{*}$ operation, i.e., for all $\sigma\in\Pi_{\text{cuts}}^{n}$, we have $\sigma^{*}\in\Pi_{\text{cuts}}^{n}$. Indeed, $(k(k+1)\cdots n12\cdots (k-1))^{*} = (n+2-k)\cdots (n-1)n12\cdots(n+1-k)\in \Pi_{\text{cuts}}^{n}$, for each $k\in [2,n]$ and $n\geq 2$.

\begin{theorem}
\label{th:costStar}
For any permutation $\pi$, $\cost(\pi) = \cost(\pi^{*})$. 
\end{theorem}
\begin{proof}
We will show that $\cost(\pi^{*})\leq \cost(\pi)$. The equality follows because $(\pi^{*})^{*} = \pi$, which will imply that $\cost(\pi)\leq \cost(\pi^{*})$.
Let $\pi = \pi_{1}\cdots\pi_{n}$. Consider an arbitrary iteration $\textsf{itr}$ over $\pi$, consisting of $m$ cuts associated with the permutations $\sigma_{1},\dots ,\sigma_{m}$, respectively. Let the cut $\sigma_{k}$ be applied over the segment $[i_{k},j_{k}]$ in $\pi$, for $k\in [m]$. Denote the output permutation after the iteration $\textsf{itr}$ with $\pi'$. Consider an iteration $\textsf{itr}^{*}$ over $\pi^{*}$, corresponding to $\textsf{itr}$, that also consists of $m$ cuts, given by the permutations $\sigma_{1}^{*},\dots ,\sigma_{m}^{*}$ which are applied over the segments $\pi^{*}_{n+1 - j_{k}}\cdots\pi^{*}_{n+1 - i_{k}}$, for $k\in [m]$. If $(\pi^{*})'$ is the output permutation after applying $\textsf{itr}^{*}$, then we claim that $(\pi^{*})'=(\pi')^{*}$. This implies that for any sequence of iterations $\textsf{itr}_{1},\textsf{itr}_{2},\ldots ,\textsf{itr}_{r}$ that sorts $\pi$, one would have a corresponding sequence of iterations $\textsf{itr}^{*}_{1},\textsf{itr}^{*}_{2},\ldots ,\textsf{itr}^{*}_{r}$ that sorts $\pi^{*}$, since $id_{n}^{*} = id_{n}$:
\begin{figure}[h!]
            \centering
\begin{tikzpicture}[scale=1]
        \node[align=right] at (0,1) {$\pi\xrightarrow{\textsf{itr}_{1}}\pi'\xrightarrow{\textsf{itr}_{2}}\cdots \xrightarrow{\textsf{itr}_{r}}id_{n}$};
        \node[align=right] at (0,0) {$\pi^{*}\xrightarrow{(\textsf{itr}_{1})^{*}}(\pi')^{*}\xrightarrow{(\textsf{itr}_{2})^{*}}\cdots \xrightarrow{(\textsf{itr}_{r})^{*}}(id_{n})^{*} = id_{n}$}; 
    \end{tikzpicture}
\end{figure}

Here is a concrete example of a single step for $\pi = 526314$ and an iteration $\textsf{itr}$ consisting of two cuts associated with the permutations $231$ and $21$ applied over the segments $\pi_{1}\pi_{2}\pi_{3} = 526$ and $\pi_{5}\pi_{6} = 14$, respectively. The corresponding iteration $\textsf{itr}^{*}$ over $\pi^{*} = 364152$ consists of the two cuts associated with the permutations $312 = 231^{*}$ and $21 = 21^{*}$ applied over the segments $\pi^{*}_{7-3}\pi^{*}_{7-2}\pi^{*}_{7-1} = \pi^{*}_{4}\pi^{*}_{5}\pi^{*}_{6} = 152$ and $\pi^{*}_{7-6}\pi^{*}_{7-5} = \pi^{*}_{1}\pi^{*}_{2} = 36$.
\begin{figure}[H]
            \centering
\begin{tikzpicture}[scale=1]
        \node[align=right] at (0,1) {$\pi = 526314\xrightarrow{\textsf{itr}}265341 = \pi'$};
        \node[align=right] at (0,0) {$\pi^{*} = 364152\xrightarrow{\textsf{itr}^{*}}634215=(\pi^{*})' = (\pi')^{*}$}; 
    \end{tikzpicture}
\end{figure}
 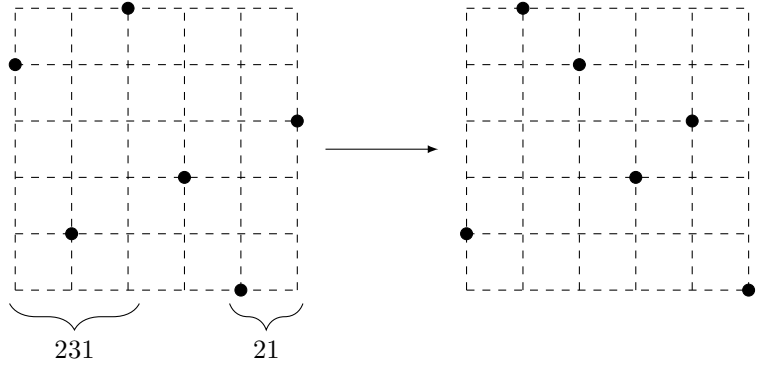
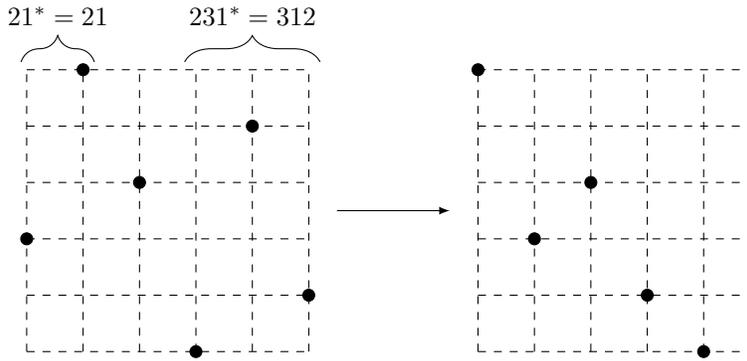
\begin{figure}[ht!]
            \centering
\begin{subfigure}[t]{\textwidth}
\captionsetup{width=\textwidth}
\centering
   \begin{tikzpicture}[scale=0.75]
   \filldraw [black] (1,5) circle (3pt);
   \filldraw [black] (2,2) circle (3pt);
   \filldraw [black] (3,6) circle (3pt);
   \filldraw [black] (4,3) circle (3pt);
   \filldraw [black] (5,1) circle (3pt);
   \filldraw [black] (6,4) circle (3pt);
 
 \draw[dashed] (1,6) -- (6,6);
 \draw[dashed] (1,5) -- (6,5);
 \draw[dashed] (1,4) -- (6,4);
 \draw[dashed] (1,3) -- (6,3);
 \draw[dashed] (1,2) -- (6,2);
 \draw[dashed] (1,1) -- (6,1);
 
 \draw[dashed] (6,1) -- (6,6);
 \draw[dashed] (5,1) -- (5,6);
 \draw[dashed] (4,1) -- (4,6);
 \draw[dashed] (3,1) -- (3,6);
 \draw[dashed] (2,1) -- (2,6);
 \draw[dashed] (1,1) -- (1,6);
 
 \draw [decorate,decoration={brace,mirror,amplitude=10pt,mirror,raise=4pt},yshift=0pt]
(0.9,0.95) -- (3.2,0.95) node [black,midway,yshift=-0.75cm] {\footnotesize $231$};

\draw [decorate,decoration={brace,mirror,amplitude=10pt,mirror,raise=4pt},yshift=0pt]
(4.8,0.95) -- (6.1,0.95) node [black,midway,yshift=-0.75cm] {\footnotesize $21$};

\draw [->] (6.5,3.5) -- (8.5,3.5);

   \filldraw [black] (9,2) circle (3pt);
   \filldraw [black] (10,6) circle (3pt);
   \filldraw [black] (11,5) circle (3pt);
   \filldraw [black] (12,3) circle (3pt);
   \filldraw [black] (13,4) circle (3pt);
   \filldraw [black] (14,1) circle (3pt);
 
 \draw[dashed] (9,6) -- (14,6);
 \draw[dashed] (9,5) -- (14,5);
 \draw[dashed] (9,4) -- (14,4);
 \draw[dashed] (9,3) -- (14,3);
 \draw[dashed] (9,2) -- (14,2);
 \draw[dashed] (9,1) -- (14,1);
 
 \draw[dashed] (9,1) -- (9,6);
 \draw[dashed] (10,1) -- (10,6);
 \draw[dashed] (11,1) -- (11,6);
 \draw[dashed] (12,1) -- (12,6);
 \draw[dashed] (13,1) -- (13,6);
 \draw[dashed] (14,1) -- (14,6);
    \end{tikzpicture}
     \captionof{figure}{The graphs of $\pi = 526314$, $\pi' = 265341$ and the action of $\textsf{itr}$}
\end{subfigure}%
\vspace{3mm}
\begin{subfigure}[t]{\textwidth}
\captionsetup{width=\textwidth}
\centering
   \begin{tikzpicture}[scale=0.75]
   \filldraw [black] (1,3) circle (3pt);
   \filldraw [black] (2,6) circle (3pt);
   \filldraw [black] (3,4) circle (3pt);
   \filldraw [black] (4,1) circle (3pt);
   \filldraw [black] (5,5) circle (3pt);
   \filldraw [black] (6,2) circle (3pt);
 
 \draw[dashed] (1,6) -- (6,6);
 \draw[dashed] (1,5) -- (6,5);
 \draw[dashed] (1,4) -- (6,4);
 \draw[dashed] (1,3) -- (6,3);
 \draw[dashed] (1,2) -- (6,2);
 \draw[dashed] (1,1) -- (6,1);
 
 \draw[dashed] (6,1) -- (6,6);
 \draw[dashed] (5,1) -- (5,6);
 \draw[dashed] (4,1) -- (4,6);
 \draw[dashed] (3,1) -- (3,6);
 \draw[dashed] (2,1) -- (2,6);
 \draw[dashed] (1,1) -- (1,6);
 
 \draw [decorate,decoration={brace,amplitude=10pt,raise=4pt},yshift=0pt]
(3.8,5.95) -- (6.2,5.95) node [black,midway,yshift=0.75cm] {\footnotesize $231^{*} = 312$};

\draw [decorate,decoration={brace,amplitude=10pt,raise=4pt},yshift=0pt]
(0.9,5.95) -- (2.2,5.95) node [black,midway,yshift=0.75cm] {\footnotesize $21^{*}= 21$};

\draw [->] (6.5,3.5) -- (8.5,3.5);

   \filldraw [black] (9,6) circle (3pt);
   \filldraw [black] (10,3) circle (3pt);
   \filldraw [black] (11,4) circle (3pt);
   \filldraw [black] (12,2) circle (3pt);
   \filldraw [black] (13,1) circle (3pt);
   \filldraw [black] (14,5) circle (3pt);
 
 \draw[dashed] (9,6) -- (14,6);
 \draw[dashed] (9,5) -- (14,5);
 \draw[dashed] (9,4) -- (14,4);
 \draw[dashed] (9,3) -- (14,3);
 \draw[dashed] (9,2) -- (14,2);
 \draw[dashed] (9,1) -- (14,1);
 
 \draw[dashed] (9,1) -- (9,6);
 \draw[dashed] (10,1) -- (10,6);
 \draw[dashed] (11,1) -- (11,6);
 \draw[dashed] (12,1) -- (12,6);
 \draw[dashed] (13,1) -- (13,6);
 \draw[dashed] (14,1) -- (14,6);
    
    \end{tikzpicture}
  \caption{The graphs of $\pi^{*} = 364152$, $(\pi^{*})' = 634215$ and the action of $\textsf{itr}^{*}$}
\end{subfigure}
\caption{Example appearing in the proof of Theorem \ref{th:costStar}. Rotate the permutation graphs on subfigure $(A)$ at $180\degree$ to obtain the permutation graphs on subfigure $(B)$.}
             \label{fig:2graphs}
        \end{figure}
$[$Proof that $(\pi')^{*} = (\pi^{*})']$ 

The graph of $\pi^{*}$ is obtained from the graph of $\pi$ by rotating at $180\degree$. In addition, the iteration $\textsf{itr}^{*}$ applies the same cuts as the iteration $\textsf{itr}$, but over the rotated graph of $\pi^{*}$ (see Figure \ref{fig:2graphs}). Therefore, the graphs of $(\pi^{*})'$ and $\pi'$ differ by a $180\degree$ rotation. Thus, if we rotate the graph of $\pi'$ by $180\degree$, we will get the same graphs, i.e., $(\pi')^{*} = (\pi^{*})'$.
\end{proof}
Note that Theorem \ref{th:costStar} holds for any shuffling method $\Sigma$ that is closed under the $^{*}$ operation.

\section{Sorting by pop shuffle queues}
\label{sec:popDev}
 One can easily see that every shuffle queue of type $(ii)$ (every pop shuffle queue) can always sort at least as many permutations as the shuffle queue of type $(i)$ for the same shuffling method. For instance, we saw, at the beginning of Section \ref{sec:cuts},  that $p_{n}(\mathbb{Q}_{\text{cuts}}^{\prime})=\mathcal{O}(d^{n}),$ where $d=2.32$. It turns out that with the pop shuffle queue for cuts, one can sort all $n!$ permutations in $S_{n}$. Below, we prove a more general statement.
\begin{theorem}
\label{th:cutsAll}
If $\Sigma$ is a shuffling method such that $(\Pi_{\Sigma}^{k})^{-1}$ contains at least one permutation ending in $j$, for every $j\in [k-1]$ and every $k\geq 2$, then
\begin{equation*}
    S_{n}(\mathbb{Q}_{\Sigma}^{\textsf{pop}}) = S_{n},
\end{equation*}
for every $n\geq 2$.
In addition, $\mathbb{Q}_{\Sigma}^{\textsf{pop}}$ can sort every permutation using a single pop operation.
\end{theorem}
\begin{proof}
We will use induction on $n$, relying on a neat observation allowing us to make the induction step. Note that $S_{2}(\mathbb{Q}_{\Sigma}^{\textsf{pop}}) = S_{2}$, since $(\Pi_{\Sigma}^{2})^{-1}$ must contain $21$ and the permutation $12$ obviously belongs to $S_{2}(\mathbb{Q}_{\Sigma}^{\textsf{pop}})$. Assume that $n>2$ and that the statement is true for all $n'<n$. Take an arbitrary permutation $\pi$ with last element $x\in [n]$ and prefix $\pi'$, i.e., $\pi = \pi'x\in S_{n}$. If $x=n$, then by the induction hypothesis, one can sort $\pi'$ and then simply push and pop $n$ to sort $\pi$. If $x\neq n$, then we know that there exists $\sigma\in (\Pi_{\Sigma}^{n})^{-1}$ ending with $x$. Take one such $\sigma$ and let $\sigma \coloneqq \sigma'x$. If we can get output $\sigma'$ on input $\pi'$ using $\mathbb{Q}_{\Sigma}^{\textsf{pop}}$ and only one pop operation, then $\pi$ would also be sortable by $\mathbb{Q}_{\Sigma}^{\textsf{pop}}$ and only one pop operation since one can get $\sigma'$ in the device, push $x$, and shuffle by applying $\sigma^{-1}$. 

However, the induction hypothesis gives us that for any input of size $n-1$, we can always get the identity as an output. In order to get output $\sigma'$ on input $\pi'$, we can relabel the elements of $\pi'$ with $1,\ldots ,n$ by looking at $\sigma'$ as the identity. Formally, since both $\pi'$ and $\sigma'$ are permutations of $[n]\setminus \{x\}$, let $\tau\in S_{n-1}$ be the permutation satisfying $\tau\pi' = \sigma'$. By the induction hypothesis, $\tau$ can be sorted by $\mathbb{Q}_{\Sigma}^{\textsf{pop}}$ using only one pop operation. Let $\textsf{itr}$ be one such iteration that sorts $\tau$ with a single pop operation. Observe that if we apply the same sequence of operations and permutations as in $\textsf{itr}$ to input $\pi'$, we will get $\sigma'$.
\begin{example}
    Consider shuffling by cuts and let $\pi = 25143$. Then, $\sigma = 45123$ is a permutation in $\Pi_{\text{cuts}}^{-1}$ that ends with $3$. We have $\pi' = 2514$, $\sigma' = 4512$. The solution to $\tau\pi' = \sigma'$ is $\tau = 4231$. Below is an iteration of $\mathbb{Q}_{\text{cuts}}^{\textsf{pop}}$ that sorts $\tau$.

\begin{equation*}
\begin{split}
    \begin{pmatrix}
4231\\
\varepsilon\\
\varepsilon\\
\end{pmatrix}
\xrightarrow{\textsf{push}}
\begin{pmatrix}
231\\
4\\
\varepsilon\\
\end{pmatrix}
\xrightarrow{\textsf{push}}
\begin{pmatrix}
31\\
42\\
\varepsilon\\
\end{pmatrix}
\xrightarrow{\textsf{push}}
\begin{pmatrix}
1\\
423\\
\varepsilon\\
\end{pmatrix} 
\xrightarrow{\substack{\textsf{shuffle}\\(cut)}}
\begin{pmatrix}
1\\
234\\
\varepsilon\\
\end{pmatrix} \\[20pt]
\xrightarrow{\textsf{push}}
\begin{pmatrix}
\varepsilon\\
2341\\
\varepsilon\\
\end{pmatrix}
\xrightarrow{\substack{\textsf{shuffle}\\(cut)}}
\begin{pmatrix}
\varepsilon\\
1234\\
\varepsilon\\
\end{pmatrix}
\xrightarrow{\substack{\textsf{pop}\\(unload)}}
\begin{pmatrix}
\varepsilon\\
\varepsilon\\
1234\\
\end{pmatrix}
\end{split}
\end{equation*}
    The same sequence of operations and permutations applied on each shuffle will give an output $\sigma' = 4512$ on input $\pi' = 2514$ :
    
    \begin{equation*}
\begin{split}
    \begin{pmatrix}
2514\\
\varepsilon\\
\varepsilon\\
\end{pmatrix}
\xrightarrow{\textsf{push}}
\begin{pmatrix}
514\\
2\\
\varepsilon\\
\end{pmatrix}
\xrightarrow{\textsf{push}}
\begin{pmatrix}
14\\
25\\
\varepsilon\\
\end{pmatrix}
\xrightarrow{\textsf{push}}
\begin{pmatrix}
4\\
251\\
\varepsilon\\
\end{pmatrix} 
\xrightarrow{\substack{\textsf{shuffle}\\(cut)}}
\begin{pmatrix}
4\\
512\\
\varepsilon\\
\end{pmatrix} \\[20pt]
\xrightarrow{\textsf{push}}
\begin{pmatrix}
\varepsilon\\
5124\\
\varepsilon\\
\end{pmatrix}
\xrightarrow{\substack{\textsf{shuffle}\\(cut)}}
\begin{pmatrix}
\varepsilon\\
4512\\
\varepsilon\\
\end{pmatrix}
\xrightarrow{\substack{\textsf{pop}\\(unload)}}
\begin{pmatrix}
\varepsilon\\
\varepsilon\\
4512\\
\end{pmatrix}
\end{split}
\end{equation*}
\end{example}
\end{proof}
Since the shuffling by $cuts$ satisfies the condition described in Theorem \ref{th:cutsAll}, we get the following corollary.
\begin{corollary}
\label{cor:cutsPop}
\begin{equation*}
    p_{n}(\mathbb{Q}_{\text{cuts}}^{\textsf{pop}}) = n!.
\end{equation*}
\end{corollary}
\subsection{Pop shuffle queues for back-front shuffling methods}
\label{sec:backFront}
In this subsection, we prove Theorem \ref{th:selfInd}, which is an analogue of Theorem \ref{th:irredPerms} for pop shuffle queues. However, Theorem \ref{th:selfInd} holds for a smaller set of shuffling methods compared to Theorem \ref{th:irredPerms}, which requires the corresponding shuffling method to have a permutation family consisting of sets of irreducible permutations. Some examples show that if we consider the same collection of shuffling methods for pop shuffle queues, we would not have a similar one-to-one correspondence as in the proof of Theorem \ref{th:irredPerms}. Nevertheless, we have such a correspondence if we constrain ourselves to shuffling methods having a stronger property which we call the \emph{back-front} property.

\begin{definition}[Back-front shuffling method]
A shuffling method $\Sigma$ is \emph{back-front} if for every $n\geq 2$, $|\Pi_{\Sigma}^{n}| = 1$, i.e., $\Pi_{\Sigma}^{n} = \{\sigma_{n}\}$ for some $\sigma_{n}\in S_{n}$ and $\sigma_{n}$ begins with $n$, i.e., the card at the back always goes at the front.
\end{definition}

One shuffling method having this property is the $\text{rev}$ method defined in Section \ref{sec:deque}. Another example is the shuffling method \emph{\text{top-bottom}} that is defined below. This method simply switches the top and bottom card.

\begin{definition} The shuffling method $top\text{-}bottom$:
\begin{equation*}
\label{eq:topBottom}
\forall n\geq 2: \Pi_{top\text{-}bottom}^{n} = \{n23\cdots (n-1)1\}.    
\end{equation*}
\end{definition}
\begin{example}
\label{ex:topBottom}
    Consider an iteration of $\mathbb{Q}_{top\text{-}bottom}^{\textsf{pop}}$ over $32415$.
\begin{equation*}
\begin{split}
    \begin{pmatrix}
$32415$\\
\varepsilon\\
\varepsilon\\
\end{pmatrix}
\xrightarrow{\textsf{push}}
\begin{pmatrix}
2415\\
3\\
\varepsilon\\
\end{pmatrix}
\xrightarrow{\textsf{push}}
\begin{pmatrix}
415\\
32\\
\varepsilon\\
\end{pmatrix}
\xrightarrow{\textsf{push}}
\begin{pmatrix}
15\\
324\\
\varepsilon\\
\end{pmatrix}
\xrightarrow{\textsf{shuffle}}
\begin{pmatrix}
15\\
423\\
\varepsilon\\
\end{pmatrix} 
\xrightarrow{\textsf{push}} 
\begin{pmatrix}
5\\
4231\\
\varepsilon\\
\end{pmatrix} \\[20pt]
\xrightarrow{\textsf{shuffle}}
\begin{pmatrix}
5\\
1234\\
\varepsilon\\
\end{pmatrix} 
\xrightarrow{\substack{\textsf{pop}\\(unload)}}
\begin{pmatrix}
5\\
\varepsilon\\
1234\\
\end{pmatrix}
\xrightarrow{\textsf{push}}
\begin{pmatrix}
\varepsilon\\
5\\
1234\\
\end{pmatrix}
\xrightarrow{\substack{\textsf{pop}\\(unload)}}
\begin{pmatrix}
\varepsilon\\
\varepsilon\\
12345\\
\end{pmatrix}\\[5pt]
\end{split}
\end{equation*}
\begingroup\vspace*{-\baselineskip}
 \label{fig:itTopBottom}
\vspace*{\baselineskip}\endgroup
\end{example}

\begin{theorem}
\label{th:selfInd}
For every back-front shuffling method $\Sigma$ and every $n\geq 2$,
\begin{equation}
p_{n}(\mathbb{Q}_{\Sigma}^{\textsf{pop}}) = F_{2n-1},    
\end{equation}
where $F_{i}$ is the $i$-th Fibonacci number with $F_{1}=F_{2} = 1$.
\end{theorem}
\begin{proof}
First, recall that we do not allow two consecutive shuffle operations when using shuffle queues. Therefore, the output after an iteration of $\mathbb{Q}_{\Sigma}^{\textsf{pop}}$ over a permutation $\pi$ is determined by the list of segments of $\pi$ that were shuffled, since $|\Pi_{\Sigma}^{m}| = 1$ for every $m\geq 2$. For instance, the list of shuffled segments for the iteration of $\mathbb{Q}_{top\text{-}bottom}^{\textsf{pop}}$ over $324165$ shown in Example \ref{ex:topBottom} is $([1,3],[1,4],[5,6])$, since exactly three shuffle operations were performed and the device contained the corresponding segment of $\pi$ before each of them, respectively. A list of shuffled segments $l$ will always be in lexicographical order, i.e., $l = ([a_{1},b_{1}],\ldots ,[a_{r},b_{r}])$, where $a_{i}\leq a_{j}$ whenever $i<j$ and $b_{u} < b_{v}$ whenever $a_{u} = a_{v}$ and $u<v$. Since we are using pop shuffle queues, if two segments overlap, they must have the same beginning. Thus when describing a list of segments with the same beginning, we will use the shorthand $[a;b_{1},b_{2},\ldots , b_{v}]$ to denote $[a,b_{1}],[a,b_{2}],\ldots ,[a,b_{v}]$ and we will call such a list a \emph{cluster}.

We are interested in the possible lists of shuffled segments when sorting a permutation with $\mathbb{Q}_{\Sigma}^{\textsf{pop}}$ or equivalently in the possible lists of clusters. Denote this set of possible lists of clusters for input of size $n$ by $LC_{n}$. Note that the set $LC_{n}$ does not depend on the shuffling method.

The idea of this proof is to show that for any $\pi\in S_{n}(\mathbb{Q}_{\Sigma}^{\textsf{pop}})$, there exists a single list of clusters in $LC_{n}$, such that any iteration over $\mathbb{Q}_{\Sigma}^{\textsf{pop}}$ corresponding to it sorts $\pi$ and vice versa - for any given list of clusters in $LC_{n}$, there exists a single $\pi\in S_{n}(\mathbb{Q}_{\Sigma}^{\textsf{pop}})$ that can be sorted by the iterations corresponding to this list, i.e., by shuffling the given clusters. This will establish a one-to-one correspondence between the sets $S_{n}(\mathbb{Q}_{\Sigma}^{\textsf{pop}})$ and $LC_{n}$. Then, we will show that $|LC_{n}|=F_{2n-1}$.

$[LC_{n}$ $\rightarrow$ $S_{n}(\mathbb{Q}_{\Sigma}^{\textsf{pop}})]$ Let $\Pi_{\Sigma}^{m} = \{\sigma_{m}\}$ and let $l\in LC_{n}$. Assume that $l=([a_{1},b_{1}],\ldots ,[a_{r},b_{r}])$. Take $id_{n}$ and apply consecutively $\sigma_{b_{r}-a_{r}}^{-1}$ over the segment $[a_{r},b_{r}]$, $\sigma_{b_{r-1}-a_{r-1}}^{-1}$ over the segment $[a_{r-1},b_{r-1}]$ and so on. After applying $\sigma_{b_{1}-a_{1}}^{-1}$ over $[a_{1},b_{1}]$, we will obtain a permutation $\pi\in S_{n}$, which can obviously be sorted by $\mathbb{Q}_{\Sigma}^{\textsf{pop}}$ by shuffling the segments of $\pi$ in the list $l$. There cannot be another such permutation $\pi'$ since $\pi'$ must be obtained from $id_{n}$ after applying consecutively $\sigma_{b_{j}-a_{j}}^{-1}$ over the segment $[a_{j},b_{j}]$, for $j = r, r-1, \ldots , 1$. There is only one permutation that can be obtained in this way and therefore for every $l\in LC_{n}$, we have only one $\pi\in S_{n}(\mathbb{Q}_{\Sigma}^{\textsf{pop}})$ that can be sorted with the iterations corresponding to $l$.

$[S_{n}(\mathbb{Q}_{\Sigma}^{\textsf{pop}})$ $\rightarrow$ $LC_{n}]$
Let $\pi\in S_{n}(\mathbb{Q}_{\Sigma}^{\textsf{pop}})$ and let us assume that $\pi$ can be sorted by two different iterations $it_{1}$ and $it_{2}$ over $\mathbb{Q}_{\Sigma}^{\textsf{pop}}$ corresponding to two different lists of clusters in $LC_{n}$, denoted by $l_{1}$ and $l_{2}$ with their last clusters denoted by $[a_{1};b_{11},\dots ,b_{1c_{1}}]$ and $[a_{1}';b_{11}',\dots ,b_{1d_{1}}']$, respectively. If these last clusters are the same, then before applying the shuffles in each of them, we must have the same permutation in the content of the device. Therefore, without loss of generality, we can assume that $[a_{1};b_{11},\dots ,b_{1c_{1}}] \neq [a_{1}';b_{11}',\dots ,b_{1d_{1}}']$ and that $[a_{1},b_{1c_{1}}]\neq [a_{1}',b_{1d_{1}}']$. If $b_{1c_{1}}\neq b_{1d_{1}}'$, then let $b_{1c_{1}}< b_{1d_{1}}'$, without loss of generality. In this case, we will have that the element $\pi_{b_{1d_{1}}'}$ is not moved anywhere when sorting $\pi$ by $it_{1}$ (and thus $\pi_{b_{1d_{1}}'}=b_{1d_{1}}'$). However, when sorting $\pi$ by $it_{2}$, $\pi_{b_{1d_{1}}'}$ is moved at position $a_{1}'$ since the method $\Sigma$ is back-front and this is the last time this element is moved. Note also that $a_{1}'\neq b_{1d_{1}}'$. Therefore, $it_{2}$ does not sort $\pi$, which is a contradiction.

Thus $b_{1c_{1}} = b_{1d_{1}}'$ and we must have that $a_{1} \neq a_{1}'$. Let $x\coloneqq b_{1c_{1}} = b_{1d_{1}}'$. Then, $\pi_{x}$ goes to position $a_{1}$ and $a_{1}'$, when we sort $\pi$ with $it_{1}$ and $it_{2}$, respectively. This means that $\pi_{x}=a_{1}$ and that $\pi_{x}=a_{1}'$, but $a_{1} \neq a_{1}'$. This is a contradiction, which shows that any $\pi\in S_{n}(\mathbb{Q}_{\Sigma}^{\textsf{pop}})$ can be sorted by iterations over $\mathbb{Q}_{\Sigma}^{\textsf{pop}}$ corresponding to exactly one list of clusters in $LC_{n}$. 

$[$Finding $|LC_{n}|]$
The desired correspondence between  $S_{n}(\mathbb{Q}_{\Sigma}^{\textsf{pop}})$ and $LC_{n}$ was established. Therefore, it suffices to get the number of different possible lists of clusters, $|LC_{n}|$, in order to find $p_{n}(\mathbb{Q}_{\Sigma}^{\textsf{pop}})$. If $l=([a_{1};b_{11},\dots ,b_{1c_{1}}],\dots ,[a_{m};b_{m1},\dots ,b_{mc_{m}}])\in LC_{n}$, then $l$ is determined by the $c_{1}'+\cdots + c_{m}' \coloneqq k$ numbers in $[n]$ comprising $l$, where we have $c_{j}'\coloneqq c_{j}+1$ numbers in the $j$-th cluster and $c_{j}'\geq 2$ for each $j\in [m]$. We have $a_{1}<b_{11}< \cdots < b_{1c_{1}}< a_{2}<\cdots < a_{m} < b_{m1} < \cdots < b_{mc_{m}}$ and thus these $k$ numbers can be chosen in $\binom{n}{k}$ ways, where $k\in [2,n]$. The number of compositions $c_{1}'+\cdots +c_{m}' = k$, where each $c_{j}'\geq 2$ is $F_{k-1}$, as proved in Lemma \ref{lemma:comp2} following this proof. When $k=0$, we have the empty set of clusters. Therefore, we obtain $|LC_{n}| = 1+ \sum\limits_{k=2}^{n}\binom{n}{k}F_{k-1}$, which is shown to be equal to $F_{2n-1}$ in Lemma \ref{lemma:fiboEq} via a nice combinatorial argument. 
\end{proof}
We include proofs of the following two lemmas for the sake of completeness.
\begin{lemma}(\cite[Exercise 1.35b]{stanleyEC1})
\label{lemma:comp2}
The number of compositions
\begin{equation*}
    \alpha_{1} + \alpha_{2}+ \cdots + \alpha_{m} = k,
\end{equation*}
of an integer $k\geq 2$, where each part $\alpha_{j}\geq 2$, is given by the Fibonacci number $F_{k-1}$, where $F_{1}=F_{2}=1$.
\end{lemma}
\begin{proof}
We will use induction. One such composition exists when $k=2$ or $k=3$. Assume that the statement holds for all integers less than some $k \geq 4$. Consider the first part $\alpha_{1}$ of a composition $\alpha_{1} + \cdots + \alpha_{m} = k$ having parts $\alpha_{j}\geq 2$ for all $j\in [m]$ and some $m\geq 1$. If $\alpha_{1}=2$, then by the induction hypothesis, there are $F_{k-3}$ compositions, $\alpha_{2} + \cdots + \alpha_{m} = k-2$. If $\alpha_{1}>2$, then take $\alpha_{1}' \coloneqq \alpha_{1}-1$. It will suffice to find the number of compositions $\alpha_{1}' + \alpha_{2}+ \cdots +\alpha_{m} = k-1$. We have that $\alpha_{1}'\geq 2$, so the number of these compositions is $F_{k-2}$, by the induction hypothesis. In total, we have $F_{k-2} + F_{k-3} = F_{k-1}$ compositions. 
\end{proof}

\begin{lemma}(\cite[Chapter 1, Identity 20]{benjamin})
\label{lemma:fiboEq}
\begin{equation}
    F_{2n-1} = 1 + \sum\limits_{k=2}^{n}\binom{n}{k}F_{k-1},
\end{equation}
where $F_{j}$ denotes the $j$-th Fibonacci number and $F_{1} = F_{2} = 1$.
\begin{proof}
It is a well-known fact that $F_{2n-1}$ is the number of ways of tiling a strip of length $2n-2$ with tiles of length either $1$ (squares) or $2$ (dominoes). If $0$ of the tiles are squares, then we have a single possible tiling of the strip. It is not possible to have exactly $1$ square, since $2n-2$ is even. The cases when the number of squares is two or more remains. In these cases, the total number of tiles is not less than $2 + \frac{2n-4}{2} = n$ and it is not possible to have less than two squares among any $n$ of the tiles. Let $k$ be the number of squares among the first $n$ tiles. We have $k\geq 2$ and there are $\binom{n}{k}$ ways of arranging these first $n$ tiles. This will be a collection of $k$ squares and $n-k$ dominoes, which has total length $2n-k$. The strip that remains has length $k-2$, which means that it can be tiled in $F_{k-1}$ ways. Summation over $k$ gives the desired sum. \end{proof}
\end{lemma}

\subsection{A conjecture on Wilf-pop-equivalence}
\label{sec:conj}
Theorem \ref{th:selfInd} enumerates the number of sortable permutations for a small subset of pop shuffle queues. One can consider sorting by pop shuffle queues for various other shuffling methods common in the literature. Some of them are the \emph{In-shuffles} and the \emph{Out-shuffles}, as well as the \emph{Monge shuffles} and the \emph{Milk shuffles}. The book of Diaconis and Graham \cite[Chapter 6]{diaconisGraham} discusses these and other shuffling methods. There, the authors mention an interesting connection between the Out-shuffles and the Milk shuffles, shown to them by John Conway. The same connection exists for the In-shuffles and the Monge shuffles. In this section, we conjecture and investigate another possible connection between the last two shuffling methods that is related to their pop shuffle queues. In order to describe it, we first define the two methods formally.

The \emph{In-shuffle} method is one of the two kinds of perfect riffle shuffles that are probably the most popular shuffling methods. When using the perfect riffle shuffles, half of the deck is held in each hand with the thumbs inward, then cards are released by the thumbs so that they fall to the table interleaved perfectly, i.e., the first card is coming from one of the halves, the second from the other half and so on. The Out-shuffles leaves the original top card back on top. The In-shuffles leaves the original top card second from top. For example, a deck of eight cards numbered by $1,2,3,4,5,6,7,8$ from top to bottom, is transformed to $5,1,6,2,7,3,8,4$ after one In-shuffle. Some applications of the riffle shuffles and some of their mathematical properties are discussed in \cite{atkinson, perfSh}. The permutation family of the In-shuffle method is 
\begin{equation*}
        \forall n\geq 2: \Pi_{\text{In-sh}}^{n} = \begin{cases}
\{(k+1)1(k+2)2 \cdots (2k)k\}, \text{if $n=2k$, and}\\
\{(k+1)1(k+2)2 \cdots (2k)k(2k+1)\}, \text{if $n=2k+1$.}
\end{cases}
    \end{equation*}
    
The \emph{Monge shuffle} method is named after the eighteenth-century geometer Gaspard Monge, who worked out the basic mathematical details of these shuffles \cite{diaconisGraham}. The Monge shuffle is carried out by successively putting cards over and under. The top card is taken into the other hand, the next is placed above, the third below these two cards and so on. For example, a deck of eight cards numbered by $1,2,3,4,5,6,7,8$ from top to bottom, is transformed to $8,6,4,2,1,3,5,7$ after one Monge shuffle. The permutation family of the Monge shuffling method is
\begin{equation*}
        \forall n\geq 2: \Pi_{\text{Monge}}^{n} = \{ \cdots 642135 \cdots\}.
    \end{equation*}
\begin{definition}[Wilf-pop-equivalent shuffling methods]
The shuffling methods $\Sigma_{1}$ and $\Sigma_{2}$ are \emph{Wilf-pop-equivalent} if for each $n\geq 1$,
\begin{equation*}
p_{n}(\mathbb{Q}_{\Sigma_{1}}^{\textsf{pop}}) = p_{n}(\mathbb{Q}_{\Sigma_{2}}^{\textsf{pop}}).
\end{equation*}
\end{definition}
We make the following conjecture.

\begin{conjecture} \label{conj:main}
The In-shuffle and the Monge shuffling methods are Wilf-pop-equivalent.
\end{conjecture}

A first step that may help to establish the conjecture is the next theorem, which confirms it if one has to use a single pop operation. Let $S_{n}^{1}(\mathbb{Q}_{\Sigma}^{\textsf{pop}})$ be the set of permutations of size $n$ sortable by $\mathbb{Q}_{\Sigma}^{\textsf{pop}}$ using only one pop operation and let $p_{n}^{1}(\mathbb{Q}_{\Sigma}^{\textsf{pop}}) \coloneqq |S_{n}^{1}(\mathbb{Q}_{\Sigma}^{\textsf{pop}})|$. 

\begin{theorem}\label{th:1pop}
For every $n\geq 1$,
\begin{equation*}
    p_{n}^{1}(\mathbb{Q}_{\text{In-sh}}^{\textsf{pop}}) = p_{n}^{1}(\mathbb{Q}_{\text{Monge}}^{\textsf{pop}}).
\end{equation*}
In addition, for every $n\geq 3$, $p_{n}^{1}(\mathbb{Q}_{\text{In-sh}}^{\textsf{pop}}) = p_{n}^{1}(\mathbb{Q}_{\text{Monge}}^{\textsf{pop}}) = a_{n-2}$, where $a_{1}=2$, $a_{2}=4$ and $a_{n}=3a_{n-2}$ for $n\geq 3$ (sequence A068911 in \cite{OEIS}).
\end{theorem}
We will show separately, with the next two lemmas, that $p_{n}^{1}(\mathbb{Q}_{\text{In-sh}}^{\textsf{pop}})$ and $p_{n}^{1}(\mathbb{Q}_{\text{Monge}}^{\textsf{pop}})$ are equal to $a_{n-2}$, for all $n > 4$. This will suffice to establish Theorem \ref{th:1pop}, since one can check directly that the statement of the theorem holds for $n \leq 4$.
\begin{lemma}
 $p_{n}^{1}(\mathbb{Q}_{\text{Monge}}^{\textsf{pop}}) = a_{n-2}$, for all $n > 4$. 
\end{lemma}
\begin{proof}
Let $\Pi_{\text{Monge}}^{i} = \{\sigma_{i}\}$, for $i>1$. We will need to use permutations of the same size. Thus let
\begin{equation}
\label{eq:tau}
\tau_{i}(x) \coloneqq
\begin{cases}
\sigma_{i}(x), & \text{if }x\leq i\text{,}\\
\text{$x$}, & \text{if $x>i$}.
\end{cases}    
\end{equation}
be a permutation of size $n$, for $i\in [2,n]$. Recall that sorting $\pi\in S_{n}$ with an iteration over $\mathbb{Q}_{\text{Monge}}^{\textsf{pop}}$ having a single pop, corresponds to a cluster $[1;b_{1},b_{2},\dots ,n]$, where one performs a shuffle after pushing $b_{1},b_{2},\dots ,n$ elements, respectively. The output will be $\pi\tau_{b_{1}}\tau_{b_{2}}\cdots\tau_{n} = id_{n}$. In general, the set of the possible iterations with a single pop over $\mathbb{Q}_{\text{Monge}}^{\textsf{pop}}$ is described by the set of vectors $(\delta_{2},\dots ,\delta_{n})$, where $\delta_{i}=0$ or $1$, for each $i\in [2,n]$ and the set of possible outputs on input $\pi$ is described by $\pi\tau_{2}^{\delta_{2}}\cdots\tau_{n}^{\delta_{n}}$. Note that if $2j+1<n$ and $j\geq 1$, then $\tau_{2j} = \tau_{2j+1}$. 

Therefore, if $n=2k+1$, for a given $k$, then the possible outputs on input $\pi$ are $\pi\tau_{2}^{\delta'_{2}}\tau_{4}^{\delta'_{4}}\cdots\tau_{2k}^{\delta'_{2k}}$, where $\delta'_{2i} = 0,1$ or $2$, for each $i\in [k]$. If $\pi\in S_{n}^{1}(\mathbb{Q}_{\text{Monge}}^{\textsf{pop}})$, then $\pi$ is a solution to $\pi\tau_{2}^{\delta'_{2}}\tau_{4}^{\delta'_{4}}\cdots\tau_{2k}^{\delta'_{2k}} = id_{n}$ for some $(\delta'_{2},\dots ,\delta'_{2k})$. Thus, $p_{n}^{1}(\mathbb{Q}_{\text{Monge}}^{\textsf{pop}})$ is given by the number of different products $\tau_{2}^{\delta'_{2}}\tau_{4}^{\delta'_{4}}\cdots\tau_{2k}^{\delta'_{2k}}$. We will show that $\tau_{2}^{\delta'_{2}}\tau_{4}^{\delta'_{4}}\cdots\tau_{2k}^{\delta'_{2k}}\neq \tau_{2}^{\delta''_{2}}\tau_{4}^{\delta''_{4}}\cdots\tau_{2k}^{\delta''_{2k}}$, if $(\delta'_{2},\dots ,\delta'_{2k})\neq (\delta''_{2},\dots ,\delta''_{2k})$. This means that it suffices to count the number of different vectors $(\delta'_{2},\dots ,\delta'_{2k})$ which implies that $p_{n}^{1}(\mathbb{Q}_{\text{Monge}}^{\textsf{pop}}) = 3p_{n-2}^{1}(\mathbb{Q}_{\text{Monge}}^{\textsf{pop}})$ since $\delta'_{2k}$ has three possible values. 

Assume that $\tau_{2}^{\delta'_{2}}\tau_{4}^{\delta'_{4}}\cdots\tau_{2k}^{\delta'_{2k}}= \tau_{2}^{\delta''_{2}}\tau_{4}^{\delta''_{4}}\cdots\tau_{2k}^{\delta''_{2k}}$ for some $(\delta'_{2},\dots ,\delta'_{2k})\neq (\delta''_{2},\dots ,\delta''_{2k})$. We can further assume that $\delta'_{2k}< \delta''_{2k}$. Therefore, $\tau_{2}^{\delta'_{2}}\tau_{4}^{\delta'_{4}}\cdots\tau_{2k-2}^{\delta'_{2k-2}}= \tau_{2}^{\delta''_{2}}\tau_{4}^{\delta''_{4}}\cdots\tau_{2k-2}^{\delta''_{2k-2}}\tau_{2k}^{\delta''_{2k}-\delta'_{2k}}$, where $\delta''_{2k}-\delta'_{2k}\in \{1,2\}$. However, $\tau_{2k}(1) = 2k$ and  $\tau_{2k}^{2}(1) = 2k-1$, i.e., $\tau_{2k}^{\delta''_{2k}-\delta'_{2k}}$ moves one of the last two elements to the first position, while neither of $\tau_{2},\dots ,\tau_{2k-2}$ moves any of these two elements, which is a contradiction. If $n=2k$, we can proceed in a similar way. We would still have $p_{n}^{1}(\mathbb{Q}_{\text{Monge}}^{\textsf{pop}}) = 3p_{n-2}^{1}(\mathbb{Q}_{\text{Monge}}^{\textsf{pop}})$ since we have three times more possibilities for the vector $(\delta_{2},\dots ,\delta_{2k-2},\delta_{2k}),$ where $\delta_{2k}$ is $0$ or $1$ and $\delta_{2i}$ is $0,1$ or $2$ for $i\in [2,k-1]$, compared to $(\delta_{2},\dots ,\delta_{2k-2}),$ where $\delta_{2k-2}$ is $0$ or $1$ and $\delta_{2i}$ is $0,1$ or $2$ for $i\in [2,k-2]$.
\newline
\end{proof}
\begin{lemma}
$p_{n}^{1}(\mathbb{Q}_{\text{In-sh}}^{\textsf{pop}}) = a_{n-2}$, for all $n>4$.
\end{lemma}
\begin{proof}
Let $\Pi_{\text{In-sh}}^{i} = \{\sigma_{i}\}$, for $i>1$ and let
\begin{equation*}
\tau_{i}(x) \coloneqq
\begin{cases}
\sigma_{i}(x), & \text{if $x\leq i$,}\\
\text{$x$}, & \text{if $x>i$}.
\end{cases}  
\end{equation*}
be a permutation of size $n$, for $i\in [2,n]$. Again, we have that $\tau_{2j} = \tau_{2j+1}$, if $2j+1\leq n$. The number $p_{n}^{1}(\mathbb{Q}_{\text{In-sh}}^{\textsf{pop}})$ is given by the number of solutions of $\pi\tau_{2}^{\delta'_{2}}\tau_{4}^{\delta'_{4}}\cdots\tau_{2k}^{\delta'_{2k}} = id_{n}$ for some $(\delta'_{2},\dots ,\delta'_{2k})$, where $\delta_{2i}'=0,1$ or $2$, for $i\in [k]$. We will show, again, that $\tau_{2}^{\delta'_{2}}\tau_{4}^{\delta'_{4}}\cdots\tau_{2k}^{\delta'_{2k}}\neq \tau_{2}^{\delta''_{2}}\tau_{4}^{\delta''_{4}}\cdots\tau_{2k}^{\delta''_{2k}}$, if $(\delta'_{2},\dots ,\delta'_{2k})\neq (\delta''_{2},\dots ,\delta''_{2k})$. Assume the opposite. Assume also that $\delta'_{2k}< \delta''_{2k}$, without loss of generality. Therefore, $\tau_{2}^{\delta'_{2}}\tau_{4}^{\delta'_{4}}\cdots\tau_{2k-2}^{\delta'_{2k-2}}= \tau_{2}^{\delta''_{2}}\tau_{4}^{\delta''_{4}}\cdots\tau_{2k-2}^{\delta''_{2k-2}}\tau_{2k}^{\delta''_{2k}-\delta'_{2k}}$, where $\delta''_{2k}-\delta'_{2k}>0$. The possible values of $\delta''_{2k}-\delta'_{2k}$ are $1$ and $2$. Now, it suffices to see that $\tau_{2k}(2k-1) = 2k$ and $\tau_{2k}^{2}(2k-3) = 2k$, while neither of $\tau_{2},\dots ,\tau_{2k-2}$ moves the element $2k$. This is a contradiction, implying that $p_{n}^{1}(\mathbb{Q}_{\text{In-sh}}^{\textsf{pop}})$ is equal to the number of different vectors $(\delta'_{2},\dots ,\delta'_{2k})$. Thus $p_{n}^{1}(\mathbb{Q}_{\text{In-sh}}^{\textsf{pop}}) = 3p_{n-2}^{1}(\mathbb{Q}_{\text{In-sh}}^{\textsf{pop}})$ for both odd and even values of $n$, in the same way as for Monge shuffles.
\end{proof}

With the next two facts, we give recurrence relations for the number of permutations in $S_{n}(\mathbb{Q}_{\text{Monge}}^{\textsf{pop}})$ that end with $n$ and that do not end with $n$. We also show that we have similar inequalities for these two subsets of $S_{n}(\mathbb{Q}_{\text{In-sh}}^{\textsf{pop}})$. Let $p_{n}^{\prime}(\mathbb{Q}_{\Sigma}^{\textsf{pop}}) = |\{\pi \in S_{n}(\mathbb{Q}_{\Sigma}^{\textsf{pop}})\mid \pi_{n}=n\}|$ and let $p_{n}^{\prime\prime}(\mathbb{Q}_{\Sigma}^{\textsf{pop}}) = |\{\pi \in S_{n}(\mathbb{Q}_{\Sigma}^{\textsf{pop}})\mid \pi_{n}\neq n\}|$, where $\Sigma$ is a shuffling method. Denote the number of elements in the device $\mathbb{D}$ before the last pop operation, for an iteration $\textsf{itr}$ over $\mathbb{D}$, by $\lps(\textsf{itr})$, which stands for \emph{last pop size}. One observation we use is that if $\pi$ can be sorted by an iteration $\textsf{itr}$ over either $\mathbb{Q}_{\text{Monge}}^{\textsf{pop}}$ or $\mathbb{Q}_{\text{In-sh}}^{\textsf{pop}}$, then the last element of $\pi$ determines whether $\lps(\textsf{itr})$ is odd or even.

\begin{theorem} 
\label{th:piNnotEq}
For every $n\geq 1$,
\begin{equation}
\label{eq:Eq}
    p_{n}^{\prime}(\mathbb{Q}_{\text{Monge}}^{\textsf{pop}}) = p_{n-1}(\mathbb{Q}_{\text{Monge}}^{\textsf{pop}}) + \frac{1}{3}\sum\limits_{j=2}^{\lfloor \frac{n-1}{2}\rfloor} p_{2j+1}^{1}(\mathbb{Q}_{\text{Monge}}^{\textsf{pop}})p_{n-(2j+1)}(\mathbb{Q}_{\text{Monge}}^{\textsf{pop}})
\end{equation}
and 
\begin{equation}
\label{ineq:Eq}
    p_{n}^{\prime}(\mathbb{Q}_{\text{In-sh}}^{\textsf{pop}}) \leq p_{n-1}(\mathbb{Q}_{\text{In-sh}}^{\textsf{pop}}) + \frac{1}{3}\sum\limits_{j=2}^{\lfloor \frac{n-1}{2}\rfloor} p_{2j+1}^{1}(\mathbb{Q}_{\text{In-sh}}^{\textsf{pop}})p_{n-(2j+1)}(\mathbb{Q}_{\text{In-sh}}^{\textsf{pop}}).
\end{equation}
\end{theorem}

\begin{proof}
Let $\pi\in S_{n}(\mathbb{Q}_{\text{Monge}}^{\textsf{pop}})$, where $\pi_{n}=n$ and let $\textsf{itr}$ be an iteration sorting $\pi$ by the given device. Let $\lps(\textsf{itr})=k$. The sequence of operations for $\textsf{itr}$ ends either with $\textsf{push}$, $\textsf{pop}$ or with $\textsf{push}$, $\textsf{shuffle}$, $\textsf{pop}$. 

In the first case, the possible prefixes $\pi' = \pi_{1}\cdots \pi_{n-1}$ are exactly the permutations in $S_{n-1}(\mathbb{Q}_{\text{Monge}}^{\textsf{pop}})$, since the iteration $\textsf{itr}$ sorts $\pi'$ and conversely if we have an iteration $\textsf{itr'}$ that sorts some $\pi'\in S_{n-1}(\mathbb{Q}_{\text{Monge}}^{\textsf{pop}})$, then $\pi'n$ would be sorted by applying $\textsf{itr'}$ and then adding the operations $\textsf{push}, \textsf{pop}$ at the end. Therefore, we have $p_{n-1}(\mathbb{Q}_{\text{Monge}}^{\textsf{pop}})$ such permutations $\pi$. 

In the second case, the last shuffle must leave the element $\pi_{n}=n$ at the same position. The latter means that $\sigma_{k}(k)=k,$ where $\Pi_{\text{Monge}}^{k}=\{\sigma_{k}\}$, which is true if and only if $k$ is odd. Let $k>1$ be a fixed odd number. Let us also have $\pi'\coloneqq \red(\pi_{n-k+1}\cdots\pi_{n}$). As in the proof of Theorem \ref{th:1pop}, we must have that $\pi'$ is a solution of the equation $\pi'\tau_{2}^{\delta_{2}}\cdots\tau_{k-1}^{\delta_{k-1}}\tau_{k}^{\delta_{k}}=id_{k}$, for a binary vector $(\delta_{2},\dots ,\delta_{k})$ and where the permutations $\tau_{j}$ are defined by Equation \eqref{eq:tau} in the same proof. Recall that $\delta_{j}=1$ if and only if the iteration $\textsf{itr}$ has a shuffle operation immediately after the $j$-th element of $\pi'$ is pushed. Take one such solution $\pi'$ corresponding to the vector $(\delta_{2},\dots ,\delta_{k})$. We have a shuffle before the last pop, which means that $\delta_{k}=1$. If $\delta_{k-1}=0$, then we must also have $\pi'\tau_{2}^{\delta_{2}}\cdots\tau_{k-2}^{\delta_{k-2}}\tau_{k-1}=id_{k}$ since $k$ is odd and $\tau_{k-1}=\tau_{k}$. This means that $\pi'$ can be sorted with an iteration ending with the operations $\textsf{push}$ and $\textsf{pop}$ and thus the same holds for $\pi$. These permutations $\pi$ were already counted in the first case. Therefore, we must have $\pi'$ for which $\delta_{k-1}=\delta_{k}=1$ in order for $\pi$ to not be yet counted. The number of these permutations $\pi'$ is the same as the number of different products $\tau_{2}^{\delta_{2}}\cdots\tau_{k-2}^{\delta_{k-2}}$, which is $p_{k-2}^{1}(\mathbb{Q}_{\text{Monge}}^{\textsf{pop}}) = \frac{1}{3}p_{k}^{1}(\mathbb{Q}_{\text{Monge}}^{\textsf{pop}})$. The permutation $\pi_{1}\cdots\pi_{n-k}$ could be any of the permutations in $S_{n-k}(\mathbb{Q}_{\text{Monge}}^{\textsf{pop}})$. Therefore, summing over all odd values of $k=2j+1$, we get an inequality similar to Equation \eqref{eq:Eq}:
\begin{equation}
\label{ineq:tmp}
    p_{n}^{\prime}(\mathbb{Q}_{\text{Monge}}^{\textsf{pop}}) \leq p_{n-1}(\mathbb{Q}_{\text{Monge}}^{\textsf{pop}}) + \frac{1}{3}\sum\limits_{j=2}^{\lfloor \frac{n-1}{2}\rfloor} p_{2j+1}^{1}(\mathbb{Q}_{\text{Monge}}^{\textsf{pop}})p_{n-(2j+1)}(\mathbb{Q}_{\text{Monge}}^{\textsf{pop}})
\end{equation}
All of the steps of the proof so far are applicable to In-shuffles as well. Therefore, we have obtained Inequality \eqref{ineq:Eq}.

It remains to show that instead of Inequality \eqref{ineq:tmp} one can write an equality. This is true because of the following observation. Assume that $\pi$ can be sorted by $\mathbb{Q}_{\text{Monge}}^{\textsf{pop}}$ using two different iterations $\textsf{itr}$ and $\textsf{itr'}$ with $\lps(\textsf{itr})=2j+1$ and $\lps(\textsf{itr'})=2j'+1$, where $j\neq j'$. Assume also, that $\pi_{n}=n$ and that $\pi_{1}\cdots\pi_{n-1}\notin S_{n-1}(\mathbb{Q}_{\text{Monge}}^{\textsf{pop}})$. Then, both $\textsf{itr}$ and $\textsf{itr'}$ must have two shuffle operations after pushing the elements $\pi_{n-1}$ and $\pi_{n}$, respectively. In addition, $\tau_{2v}^{2}(1)=2v-1$ for each $v\geq 1$. Therefore, since $\textsf{itr}$ sorts $\pi$, we must have $\pi_{n-1}=n-(2j+1)+1 = n-2j$. However, since $\textsf{itr'}$ sorts $\pi$, we must also have $\pi_{n-1}=n-2j'$, which is a contradiction.
\end{proof}

\begin{theorem} 
\label{th:notEnding}
For every $n\geq 1$,
\begin{equation}
\label{eq:mongeNotEq}
    p_{n}^{\prime\prime}(\mathbb{Q}_{\text{Monge}}^{\textsf{pop}}) = \sum\limits_{j=1}^{\lfloor \frac{n}{2}\rfloor} p_{2j}^{1}(\mathbb{Q}_{\text{Monge}}^{\textsf{pop}})p_{n-2j}(\mathbb{Q}_{\text{Monge}}^{\textsf{pop}})
\end{equation}
and 
\begin{equation}
\label{ineq:inShNotEq}
    p_{n}^{\prime\prime}(\mathbb{Q}_{\text{In-sh}}^{\textsf{pop}}) \leq \sum\limits_{j=1}^{\lfloor \frac{n}{2}\rfloor} p_{2j}^{1}(\mathbb{Q}_{\text{In-sh}}^{\textsf{pop}})p_{n-2j}(\mathbb{Q}_{\text{In-sh}}^{\textsf{pop}}).
\end{equation}
\end{theorem}

Equation \eqref{eq:mongeNotEq} follows from a property of the Monge shuffle, which we formulate below.

\begin{definition}[Pop-simple shuffling method]
\label{def:popS}
The shuffling method $\Sigma$ is \emph{pop-simple} if there is no permutation $\pi\in S^{1}(\mathbb{Q}_{\Sigma}^{\textsf{pop}})$, not ending with $n$, such that $\pi = \pi'\oplus \pi''$ for some $\pi'$ and $\pi''$, where $|\pi'|\geq 2$, $|\pi''| \geq 2$ and $\pi'' \in S^{1}(\mathbb{Q}_{\Sigma}^{\textsf{pop}})$.
\end{definition}

Intuitively, if a shuffling method $\Sigma$ is pop-simple and $\sigma\in S_{n}(\mathbb{Q}_{\Sigma}^{\textsf{pop}})$ does not end with $n$, then $\lps(\textsf{itr})$ has the same value for every iteration $\textsf{itr}$ of $\mathbb{Q}_{\Sigma}^{\textsf{pop}}$ sorting $\sigma$. 

\begin{lemma}\label{lemma:popSimple2}
The Monge shuffling method is pop-simple.
\end{lemma}
\begin{proof}
Suppose that there exists $\pi\in S_{n}$, not ending with $n$, such that $\pi = \pi'\oplus \pi''$ for some $\pi'$ and $\pi''$, such that $|\pi'| \geq 2$, $|\pi''| \geq 2$, and each of $\pi$ and $\pi''$ can be sorted by $\mathbb{Q}_{\Sigma}^{\textsf{pop}}$ using a single pop operation. Every iteration sorting a permutation by $\mathbb{Q}_{\text{Monge}}^{\textsf{pop}}$ that uses a single pop operation can be written as a cluster beginning with the element $1$. Consider an arbitrary cluster $[1;b_{1},b_{2},\dots , b_{v}]$ representing an iteration that sorts $\pi$. Since $\pi$ does not end with $n$, then the last shuffle must be after we push the last element, i.e., $b_{v}=n$. In addition, we must have $\sigma_{n}\neq n$, where $\Pi_{\text{Monge}}^{n} = \{\sigma\}$ and $\sigma = \sigma_{1}\cdots \sigma_{n}$. Note that $\sigma$ either begins with $n$ (when $n$ is even) or ends with $n$ (when $n$ is odd). Therefore, $\sigma$ must begin with $n$, which means that $\pi_{n} = 1$. However, since $\pi = \pi'\oplus \pi''$, we must have that $1$ is among the first $|\pi'|$ elements of $\pi$. It cannot be the last element of $\pi'$, since $\pi''$ is non-empty, which is a contradiction.
\end{proof}
\begin{proof}[Proof of Theorem \ref{th:notEnding}]
Let $\pi\in S_{n}(\mathbb{Q}_{\text{Monge}}^{\textsf{pop}})$, where $\pi_{n}\neq n$ and let $\textsf{itr}$ be an iteration sorting $\pi$ by $\mathbb{Q}_{\text{Monge}}^{\textsf{pop}}$. Note that the sequence of operations for $\textsf{itr}$ must end with $\textsf{push}, \textsf{shuffle}, \textsf{pop}$ since $\pi_{n}$ must be moved to another position. This is also the reason that if $\lps(\textsf{itr})=k$, then $k$ must be even since all permutations of odd size associated with the Monge shuffle fix its last element. The permutation $\pi_{1}\cdots\pi_{n-k}$ could be any of the permutations in $S_{n-k}(\mathbb{Q}_{\text{Monge}}^{\textsf{pop}})$. Therefore, summing over all even values of $k=2j$, we get 
\begin{equation}
\label{eq:mongeNotEqTmp}
    p_{n}^{\prime\prime}(\mathbb{Q}_{\text{Monge}}^{\textsf{pop}}) \leq \sum\limits_{j=1}^{\lfloor \frac{n}{2}\rfloor} p_{2j}^{1}(\mathbb{Q}_{\text{Monge}}^{\textsf{pop}})p_{n-2j}(\mathbb{Q}_{\text{Monge}}^{\textsf{pop}}).
\end{equation}
All of the steps of the proof so far are applicable to In-shuffles, as well, and thus Inequality \eqref{ineq:inShNotEq} can be obtained analogously.

It remains to show that instead of Inequality \eqref{eq:mongeNotEqTmp}, one can write Equation \eqref{eq:mongeNotEq}. Assume that $\pi$ can be sorted by $\mathbb{Q}_{\text{Monge}}^{\textsf{pop}}$ using two different iterations $\textsf{itr}$ and $\textsf{itr'}$ with $\lps(\textsf{itr})=2j$ and $\lps(\textsf{itr'})=2j'$, where $j>j'$. Then, both sequences $\gamma = \red(\pi_{n-2j+1}\cdots\pi_{n})$ and $\kappa = \red(\pi_{n-2j'+1}\cdots\pi_{n})$ must be permutations of $[2j]$ and $[2j']$, respectively, and they must be sortable with a single pop. However, this would imply that the Monge shuffling method is not pop-simple, because $\gamma = \gamma'\oplus \kappa$ for $\gamma' = \red(\pi_{n-2j+1}\cdots\pi_{n-2j'})$, $|\gamma'|\geq 2$, $|\kappa|\geq 2$  and $\gamma$,$\kappa\in S^{1}(\mathbb{Q}_{\text{Monge}}^{\textsf{pop}})$. This contradicts Lemma \ref{lemma:popSimple2}.
\end{proof}

If one can replace Inequalities \eqref{ineq:Eq} and \eqref{ineq:inShNotEq} with equations, then one can obtain Conjecture \ref{conj:main} using induction and Theorem \ref{th:1pop}. Inequality \eqref{ineq:inShNotEq} can be replaced by an equation if and only if the In-shuffle method is also pop-simple. There exist permutations $\pi = \pi'\oplus \pi''$ for some $\pi'$ and $\pi''$, such that $|\pi'|, |\pi''| \geq 2$ and $\pi\in S^{1}(\mathbb{Q}_{\text{In-sh}}^{\textsf{pop}})$. For instance, if $\pi = 21 \oplus 62481357$, then $\pi\in S^{1}(\mathbb{Q}_{\text{In-sh}}^{\textsf{pop}})$. However, in this example, $\pi'' = 62481357 \notin S^{1}(\mathbb{Q}_{\text{In-sh}}^{\textsf{pop}})$. We have performed computer simulations using Theorem \ref{th:notEnding} which show that there is no such permutation $\pi\in S_{n}$ for $n<20$ and thus we have an equality in \eqref{ineq:inShNotEq} for $n<20$ . Similarly, we have checked that Inequality \eqref{ineq:Eq} is an equality for $n<20$. Therefore, we have obtained that Conjecture \ref{conj:main} holds for $n<20$.

\section{Further questions}
\label{sec:questions}
 The considered sorting devices and the obtained results raise some additional questions that we list below.
 \\
\begin{enumerate}[label=\arabic*)]
    
    \item Can we use Theorem \ref{th:deque} to make progress on the long-standing problem of finding the number of permutations sortable by a deque \cite[A182216]{OEIS}? Some results on the asymptotic of these numbers can be found in \cite{price,priceGut}.  \\
    \item Can we find shuffle queues that are equivalent to the input and the output restricted deques defined in \cite{knuth}? In general, if $T$ is a set of patterns, then for which $T$ exists a shuffle queue $\mathbb{Q}_{\Sigma}$, such that $S_{n}(\mathbb{Q}_{\Sigma}) = Av_{n}(T)$, for each $n\geq 2$?\\
    \item Find characterizations in terms of pattern avoiding classes for the set of permutations of given cost. Theorem \ref{th:cost1} gives such a characterization for the set of permutations of cost $1$.
    \\
    \item  In Section \ref{sec:motivation}, we noted that there exists a deterministic linear time algorithm that sorts all of the permutations in $S_{n}(\mathbb{Q}_{\text{cuts}}^{\prime})$. Which  are  the  shuffling  methods $\Sigma$ for which there exists such a linear procedure that sorts all of the permutations in $S_{n}(\mathbb{Q}_{\Sigma}^{\prime})$?
    \\
    \item Find characterization of the shuffling methods, whose shuffle queues without restrictions or of types $(i)$ and $(ii)$, can sort all permutations in $S_{n}$? Theorem \ref{th:cutsAll} identifies one class of such shuffling methods for shuffle queues of type $(ii)$.
    \\
    \item Suppose that the expected number of random shuffles until one obtains a sorted deck of cards, beginning with a random permutation, is greater or equal asymptotically when using the shuffling method $\Sigma_{1}$, compared to the shuffling method $\Sigma_{2}$. Then, is it always true that $p_{n}(\mathbb{Q}_{\Sigma_{1}}^{\prime}) \geq p_{n}(\mathbb{Q}_{\Sigma_{2}}^{\prime})$ asymptotically?
    
\end{enumerate}
One may also, of course, consider sorting by shuffle queues for other popular shuffling methods (see \cite[Chapter 6]{diaconisGraham}), in order to find more connections with other combinatorial objects. 

\section*{Acknowledgement}
I am grateful to my advisor Bridget Tenner for the useful comments and ideas.

\end{document}